\documentclass[11pt,a4paper,english]{article}
\usepackage[toc]{appendix}
\usepackage{a4wide}
\usepackage{authblk}
\usepackage[T1]{fontenc}
\usepackage{lmodern}
\usepackage{amsmath,amssymb,amsthm,mathtools,bm}
\usepackage{geometry,microtype,enumitem}
\usepackage{xcolor}
\usepackage{hyperref}
\usepackage{booktabs}
\usepackage[
backend=biber,          %
style=numeric-comp,     %
sorting=none,           %
maxnames=3,             %
url=false,              %
eprint=false,           %
isbn=false,             %
doi=true                %
]{biblatex}
\hypersetup{colorlinks=true,linkcolor=blue!50!black,citecolor=blue!50!black,urlcolor=blue!50!black,pdftitle={Geomagnetic Q-Response Inversion as a Calderon-Type Inverse Boundary Value Problem},pdfauthor={Xinpeng Huang; Zhengyong Ren; Jingtian Tang}}
\usepackage{aliascnt}
\usepackage[capitalise]{cleveref}
\numberwithin{equation}{section}
\newtheorem{theorem}{Theorem}[section]
\newaliascnt{lemma}{theorem}
\newtheorem{lemma}[lemma]{Lemma}
\aliascntresetthe{lemma}
\newaliascnt{corollary}{theorem}
\newtheorem{corollary}[corollary]{Corollary}
\aliascntresetthe{corollary}
\newaliascnt{proposition}{theorem}
\newtheorem{proposition}[proposition]{Proposition}
\aliascntresetthe{proposition}
\theoremstyle{remark}
\newaliascnt{remark}{theorem}
\newtheorem{remark}[remark]{Remark}
\aliascntresetthe{remark}
\newcommand{\R}{\mathbb{R}}
\newcommand{\B}{\mathbb{B}_R^3}
\newcommand{\curl}{\operatorname{curl}}
\newcommand{\divg}{\operatorname{div}}

\newcommand{\Span}{\operatorname{span}}

\newcommand{\dd}{\,\mathrm d}
\newcommand{\eg}{\mathbf e}
\newcommand{\bfE}{\mathbf E}
\newcommand{\bfH}{\mathbf H}
\newcommand{\bfJ}{\mathbf J}
\newcommand{\bfA}{\mathbf A}
\newcommand{\bfP}{\mathbf p}
\newcommand{\xx}{\mathbf x}
\newcommand{\ephi}{\mathbf e_\varphi}
\newcommand{\norm}[1]{\lVert #1\rVert}

\newcommand{\inner}[2]{\langle #1,#2\rangle_{\mathrm{bil}}}
\newcommand{\sesinner}[2]{\langle #1,#2\rangle}
\title{\textbf{Geomagnetic $Q$-Response Inversion as a\\Calder\'on-Type Inverse Boundary Value Problem}}
\author{\begin{tabular}{c}
Xinpeng Huang\thanks{xinpeng.huang@csu.edu.cn} \qquad Zhengyong Ren\thanks{renzhengyong@csu.edu.cn} \qquad Jingtian Tang\thanks{jttang@csu.edu.cn}\\
\small School of Geosciences and Info-Physics, Central South University,\\
\small Lushannan Road 932, Changsha, China
\end{tabular}}
\date{}
\begin{document}
\maketitle
\begin{abstract}
Whether electrical conductivity can be uniquely determined from magnetic measurements alone is a fundamental question in geomagnetic induction sounding. For a spherical conductor in the magnetoquasistatic model, we establish an equivalence between the complete Q-response and a restricted magnetic boundary operator, connecting geomagnetic induction to Calder\'on-type inverse boundary value theory.

For axisymmetric conductivities, we derive exact representations of the zero-frequency zonal response derivatives as conductivity-weighted spatial moments. These formulas give an explicit interpretation of individual response derivatives, and completeness of the moment family yields conductivity uniqueness. For layered axisymmetric models with known shell interfaces and piecewise-constant radial dependence, we further establish uniqueness at a single nonzero frequency. Finite-dimensional subclasses also admit uniqueness and stability estimates from finitely many response coefficients at sufficiently low nonzero frequencies.

These results identify complementary mechanisms for conductivity recovery from magnetic-only observations, while uniqueness for general three-dimensional scalar conductivities from ideal magnetic-only Q-response data remains unresolved.
\end{abstract}

\section{Introduction}\label{sec:introduction}
Electrical conductivity imaging provides information about the interiors of Earth and other planetary bodies and has applications in biomedical imaging and nondestructive \mbox{evaluation~\cite{Mittelholz2026,Uhlmann2009,SoleimaniLionheart2006}}. Its mathematical formulation depends on the electromagnetic observables. In electrical impedance tomography, boundary voltages and currents lead to the classical Calder\'on inverse conductivity problem \cite{Uhlmann2009,SylvesterUhlmann1987}; magnetotellurics exploits relations between electric and magnetic fields \cite{ChaveJones2012}; and Maxwell inverse boundary value theory studies material parameters through relations between tangential electric and magnetic traces \cite{Ola1993,Kenig2011}. The common identifiability question is whether the available observations determine the interior conductivity.

Several induction methods instead rely on magnetic measurements alone, which can be recorded without electrodes or direct electrical contact with the conductor. Examples include the airborne Z-axis Tipper Electromagnetic (ZTEM) method \cite{HolthamOldenburg2010} and, more directly relevant here, global geomagnetic induction sounding. Long-term ground-observatory records and the near-global coverage of satellite magnetic measurements provide complementary information about the Earth's deep electrical structure \cite{Puthe2014,Kuvshinov2021}. Related induction principles are used to investigate other planetary bodies \cite{Mittelholz2026}. These applications motivate the question addressed here: can interior conductivity be uniquely determined from magnetic-only observations?

The answer depends on the specific magnetic response operator. Under the standard vertically incident plane-wave approximation, magnetic tipper responses vanish for a horizontally layered one-dimensional Earth and therefore cannot distinguish among such conductivity profiles \cite{HolthamOldenburg2010}. To investigate what information magnetic observations can contain, we consider an idealized experiment in which the complete magnetic excitation--response relation is available.

The matrix $Q$-response framework of global geomagnetic induction sounding provides a concrete setting for this question. Within the magnetoquasistatic approximation, the incident and induced magnetic fields in a current-free spherical shell admit spherical-harmonic potential expansions. Writing $\mathbf f$ and $\widetilde{\mathbf f}$ for the incident and induced potential coefficient sequences, respectively, we have
\begin{equation}\label{eq:introQ}
 \widetilde{\mathbf f}=\mathcal{Q}_\sigma(\omega)\mathbf f,
\end{equation}
where $\sigma$ denotes conductivity and $\omega$ is angular frequency. For radially symmetric conductivity, this operator is diagonal. Its scalar response functions are related to the classical $C$-response, and their frequency dependence connects the radial inverse problem with one-dimensional inverse spectral theory \cite{Bailey1970,Weidelt1972,Weidelt2005}. Lateral conductivity variations couple spherical-harmonic modes, and matrix $Q$-responses have been used in numerical three-dimensional mantle-conductivity inversions \cite{Puthe2014,Kuvshinov2021}. Knowing the response to every admissible incident mode specifies the magnetic input--output relation within this class of spherical induction experiments. The complete operator therefore provides an ideal benchmark for conductivity identifiability; estimating it from observations is a separate problem involving source diversity, spatial coverage, source separation, and measurement error.

Within the magnetoquasistatic model, consider a spherical conductor $\B\subset\R^3$ with no prescribed internal current sources, a real, frequency-independent scalar conductivity bounded above and below by positive constants, and the same known constant permeability inside and outside. At a fixed nonzero frequency $\omega$, we ask whether equality of the complete magnetic response operators implies equality of the conductivities:
\begin{equation}\label{eq:main_question}
 \mathcal{Q}_{\sigma_1}(\omega)=\mathcal{Q}_{\sigma_2}(\omega)
 \quad\stackrel{?}{\Longrightarrow}\quad \sigma_1=\sigma_2\ \text{a.e. in }\B.
\end{equation}
We also consider equality of the complete response operators on a nonempty open frequency interval, or equality of their derivatives at zero frequency. To our knowledge, uniqueness for general three-dimensional scalar conductivities remains unresolved in each of these data settings. Either a uniqueness theorem or a counterexample would clarify the information limits of this magnetic-only experiment.

We approach the problem through Calder\'on-type inverse boundary value theory. A key distinction from standard full-data Maxwell inverse boundary problems is the available boundary information. In these experiments, the total tangential magnetic field is constrained to be a surface gradient. At nonzero frequency, the magnetic response directly determines the surface curl of the tangential electric field, but does not directly supply its surface-gradient component. The resulting restricted boundary data do not directly meet the hypotheses of standard full-data Maxwell uniqueness theorems. We formulate the corresponding restricted magnetic boundary operator, establish its equivalence with the complete $Q$-response operator, and derive a bilinear integral identity relating response differences to conductivity differences.

For conductivities axisymmetric about a known axis, we obtain an explicit characterization of the first-order low-frequency zonal ($m=0$) response. These models allow both radial and latitudinal variations while retaining the three-dimensional magnetoquasistatic Maxwell system. The zero-frequency slopes admit exact representations as conductivity-weighted spatial moments, without a small-contrast assumption. The zero-frequency slopes of off-diagonal cross-degree coefficients supply mixed moments of the meridional structure, including moments sensitive to north--south asymmetry. This characterization gives an analytical interpretation of individual response slopes and complements the use of $Q$-responses in numerical geomagnetic inversion. Completeness of the full moment family, including both diagonal and off-diagonal slopes, implies that the complete zonal slope matrix uniquely determines any axisymmetric conductivity satisfying the stated positive upper and lower bounds, up to equality almost everywhere. Analyticity extends this conclusion to exact complete zonal response data on any nonempty open interval of physical frequencies.

A different mechanism yields uniqueness from the complete zonal response at any prescribed nonzero frequency for a narrower class of axisymmetric shell models, with known radial interfaces and conductivity independent of radius and continuous in angle within each shell. The proof combines a five-dimensional scalar reduction, high-degree boundary asymptotics, and layer stripping. For finite-dimensional shell models whose conductivity in each shell is a polynomial of bounded degree in the cosine of colatitude, the low-frequency moment mechanism additionally yields uniqueness and a model-class-dependent Lipschitz stability estimate from a finite zonal response submatrix at any sufficiently low nonzero frequency. These results establish complementary forms of identifiability within specified axisymmetric classes.

Section~\ref{sec:model} defines the magnetic data and forward model and recalls the radial benchmark. Section~\ref{sec:calderon} develops the restricted boundary formulation and integral identity. Section~\ref{sec:axisym} proves the axisymmetric results. Section~\ref{sec:discussion} discusses their geophysical interpretation and practical limitations, and Section~\ref{sec:outlook} identifies open questions concerning uniqueness for general three-dimensional conductivities, general domains, and multi-frequency inversion. The appendices give the radial inverse spectral formulation and technical details of the single-frequency argument.

\section{Spherical magnetic data and the \texorpdfstring{$Q$}{Q}-response matrix}\label{sec:model}
\subsection{Incident and induced fields on an observation sphere}
We briefly recall the spherical-harmonic representation used in global geomagnetic induction sounding. Magnetometers record the superposition of incident and induced magnetic fields. Write $\mathbb{B}_a^n=\{x\in\mathbb{R}^n:|x|<a\}$ and $\mathbb{S}^{n-1}=\partial\mathbb{B}_1^n$. Let $\B$ be the conducting ball and take $\Gamma=\partial\B$ as the reference sphere for the inverse problem. Assume that a spherical shell immediately outside $\Gamma$ is insulating and current-free. We formulate the inverse problem using magnetic data on $\Gamma$. Exact full-vector data on a complete larger concentric sphere can be referred to $\Gamma$ by harmonic continuation, provided that the intervening annulus is current-free. This correspondence does not imply stable downward continuation of noisy data.    %
In the magnetoquasistatic approximation, $\bfH=-\nabla V$ in this shell. The Gauss decomposition \cite{Glassmeier2014} gives
\begin{equation}\label{eq:gauss}
\nabla V(r,\widehat x)=\sum_{\ell\geq1}\sum_{m=-\ell}^\ell
 \nabla\left\{\left[f_{\ell m}(r/R)^\ell+\tilde{f}_{\ell m}(R/r)^{\ell+1}\right]Y_{\ell m}(\widehat x)\right\}.
\end{equation}
Here $\widehat x=x/r\in\mathbb{S}^2$, and $f_{\ell m}$ and $\tilde f_{\ell m}$ are the Gauss coefficients of the external and internal potential contributions, respectively, \emph{relative to the reference sphere $\Gamma$}. The constant incident mode produces no magnetic field, while the degree-zero decaying mode is excluded by the absence of magnetic monopoles. We use the standard spherical harmonics $Y_{\ell m}$ orthonormal on the unit sphere $\mathbb{S}^2$, with real zonal elements $Y_{\ell0}=\sqrt{(2\ell+1)/(4\pi)}P_\ell$.
On $\Gamma$, let $\nu=x/R$ be the outward unit normal and $\bfH_T=\bfH-(\bfH\cdot\nu)\nu$ the tangential magnetic field. Write $H_{\nu,\ell m}$ for the coefficient of $\bfH\cdot\nu$ in $Y_{\ell m}$ and $b_{\ell m}$ for the coefficient of $\bfH_T$ in the basis $-\nabla_{\mathbb{S}^2}Y_{\ell m}/R$. Differentiating~\eqref{eq:gauss} gives
\begin{equation}\label{eq:sep1}
 b_{\ell m}=f_{\ell m}+\tilde{f}_{\ell m},\qquad
 RH_{\nu,\ell m}=-\ell f_{\ell m}+(\ell+1)\tilde{f}_{\ell m}.
\end{equation}
Consequently,
\begin{equation}\label{eq:sep2}
 f_{\ell m}=\frac{(\ell+1)b_{\ell m}-RH_{\nu,\ell m}}{2\ell+1},\qquad
 \tilde{f}_{\ell m}=\frac{\ell b_{\ell m}+RH_{\nu,\ell m}}{2\ell+1}.
\end{equation}
For ideal full-vector data on a complete sphere, the regular and decaying potential-field vector harmonics are mutually orthogonal in the spherical $L^2$ inner product, allowing their coefficients to be separated by projection.
Throughout the idealized $Q$-response formulation below, we assume that the physical inducing field has been identified and contributions from other internal magnetic sources have been removed, so that $f$ and $\tilde f$ represent the incident and induction-related contributions for the conductor under study. On a regional spherical patch, the global orthogonality used above is no longer directly available; the resulting source-separation limitations are discussed in \cite{Huang2026}. Potential-field decompositions on more general Lipschitz surfaces are studied in \cite{Baratchart2021}.

\subsection{Magnetoquasistatic model and frequency conventions}\label{subsec:forward}
In the following, we assume the conducting region is $\B\subset\R^3$, the magnetic measurement region is insulating, and the same known constant permeability $\mu>0$ applies inside and outside. The scalar conductivity is real and frequency independent, with
\begin{equation}\label{eq:condbounds}
 0<\sigma_*\leq\sigma(x)\leq\sigma^*<\infty
 \qquad\text{for a.e. }x\in\B.
\end{equation}
We neglect displacement currents and assume there are no prescribed internal current sources. For time-harmonic fields, we use the time convention $e^{-i\omega t}$, with real angular frequency $\omega$; physical fields are the real parts of their complex phasor representations.

For the causal time-domain interpretation, assume zero initial fields and use the Laplace transform
\begin{equation}\label{eq:laplace-transform}
 \widehat h(s)=\int_0^\infty e^{-st}h(t)\dd t,
 \qquad \Re s>0,
\end{equation}
whenever the integral converges. In this setting time differentiation transforms to multiplication by $s$, without an initial-field term. The transformed Maxwell equations inside the conductor are
\begin{equation}\label{eq:maxwell}
 \curl\bfE=-s\mu\bfH,\qquad
 \curl\bfH=\bfJ=\sigma\bfE,\qquad
 \divg\bfH=0.
\end{equation}
The same equations describe harmonic phasors when $s=-i\omega$. In the absence of surface sources, the tangential electric and magnetic fields and the normal component of the magnetic flux density are continuous across the interface, while the normal conduction-current trace vanishes.

\subsection{The complex magnetic \texorpdfstring{$Q$}{Q}-response operator}\label{subsec:complex-frequency}
The geomagnetic $Q$-matrix describes the frequency-domain relation between incident and induced magnetic-potential coefficients. Its entries are generally complex-valued at real angular frequencies and can be estimated from geomagnetic time series or computed using frequency-domain forward solvers \cite{Puthe2014,Kuvshinov2021}. For fixed conductivity and frequency, the induced phasor coefficients depend linearly on the incident coefficients in \eqref{eq:gauss}:
\begin{equation}\label{eq:Qdef}
 \tilde f_{kn}(\omega)
 =\sum_{\ell,m}Q_{\sigma;kn,\ell m}^{\mathrm{phys}}(\omega)f_{\ell m}(\omega).
\end{equation}
The column indexed by $(\ell,m)$ gives the induced coefficients generated by a unit incident mode, and the row indexed by $(k,n)$ gives the contributions to that output mode. %

We use a Laplace parameterization to describe the same response analytically. Let $\mathbf F(s)$ and $\widetilde{\mathbf F}(s)$ be the temporal Laplace transforms of the incident and induced coefficient sequences. The zero-initial-state problem in Section~\ref{subsec:forward} defines
\begin{equation}\label{eq:laplace-response}
 \widetilde{\mathbf F}(s)=\mathcal Q_\sigma(s)\mathbf F(s),
 \qquad \Re s>0.
\end{equation}
Here $\mathcal Q_\sigma$ is the complete magnetic response operator and $Q_{\sigma;kn,\ell m}$ denotes a scalar matrix entry. The physical-frequency operator is
\begin{equation}\label{eq:physical-laplace}
 \mathcal Q_\sigma^{\mathrm{phys}}(\omega)
 =\lim_{\varepsilon\downarrow0}\mathcal Q_\sigma(\varepsilon-i\omega)
 =\mathcal Q_\sigma(-i\omega).
\end{equation}
The resolvent construction in Section~\ref{subsec:zero_cont} establishes analyticity in a neighborhood of the closed right half-plane, so this limit is regular. The factor $e^{-\varepsilon t}$ supplies temporal damping before the harmonic limit is taken. We suppress the superscript $\mathrm{phys}$ when the argument is explicitly $\omega$; an argument $s$ denotes the Laplace parameterization. This convention also applies to scalar entries and finite blocks.

For each fixed input--output spherical-harmonic index pair, exact complex values on any nonempty open interval of physical frequencies uniquely determine the corresponding analytic entry. Thus the physical-frequency data and the analytic response family carry the same uniqueness information for the same set of spatial indices; the precise statement follows in Lemma~\ref{lem:frequency-equivalence}. %

In this finite-conductivity model, $\mathcal Q_\sigma(0)=0$. The moments derived below describe the first derivative, with
\begin{equation}\label{eq:slope-conventions}
 \left.\partial_\omega \mathcal{Q}_{\sigma}^{\mathrm{phys}}(\omega)\right|_{\omega=0}
 = \left. -i\,\partial_s \mathcal{Q}_\sigma(s)  \right|_{s=0}
\end{equation}
where the derivatives are understood entrywise. We also write $\rho=(x_1^2+x_2^2)^{1/2}$ and $z=x_3$. For conductivity $\sigma=\sigma(\rho,z)$ symmetric about the known $z$-axis, rotation symmetry preserves the azimuthal order $m$ while allowing coupling between degrees. The block with $m=0$ at both input and output is denoted by $\mathcal Q_\sigma^{00}$, and its finite submatrix of degrees $1,\ldots,N$ by $Q_{\sigma,N}^{00}$. The boundary formulation in Section~\ref{sec:calderon} does not assume axisymmetry.

The Gauss decomposition permits separation of internal and external magnetic-potential contributions from ideal full-vector data on a complete sphere, but this separation does not by itself determine the $Q$-response matrix. Even after the physical incident and induced fields have been identified, estimation of its elements depends on the spatial diversity and temporal frequency content of the incident field \cite{Puthe2014,Puthe2015}. In practice, natural sources may not provide sufficiently independent excitations of the relevant spatial modes over the frequency range of interest. The complete $Q$-response therefore serves as an idealized benchmark, assuming knowledge of the response to every admissible incident mode, with the required frequency information specified separately for each uniqueness problem. Our focus is whether the prescribed response data uniquely determine conductivity; estimating those data from finite geomagnetic observations is a separate problem.

\subsection{Uniqueness for radial conductivity}\label{subsec:radial_benchmark}
Before considering general three-dimensional conductivities, we recall the radially symmetric case as a classical benchmark for uniqueness from magnetic induction data. If $\sigma=\sigma(r)$, rotational symmetry makes the response diagonal in spherical harmonics, with the same scalar response for all orders of a given degree:
\begin{equation}\label{eq:radialdiagonal}
 Q_{\sigma;kn,\ell m}(s)=q_\ell(s)\,\delta_{k\ell}\delta_{nm}.
\end{equation}
For any fixed degree $\ell\geq1$, the complete frequency dependence of the scalar response $q_\ell(s)$ determines the radial conductivity through a one-dimensional inverse spectral problem. In the formulation recalled in Appendix~\ref{app:radial}, the response is equivalent to the Weyl function of a positive-density string, whose inverse spectral correspondence determines $\sigma(r)$ almost everywhere. This is a uniqueness result based on the full frequency-dependent response of one spatial mode, not on a single frequency sample. Hooshyar and Razavy \cite{HooshyarRazavy1982} considered a complementary fixed-frequency setting in which responses at different spherical-harmonic degrees are used in a finite-difference and Thiele continued-fraction reconstruction. Their reconstruction is formulated for a discretized radial model. The continuous model classes and data assumptions considered below are different.

The situation changes when lateral conductivity variations are allowed: different spherical-harmonic modes become coupled, and the scalar inverse spectral reduction no longer applies directly. The corresponding question is whether the full matrix response retains enough information to determine a non-radial conductivity. A complementary radial mechanism, based instead on zero-frequency slopes across spatial degrees, will be recalled in Section~\ref{subsec:slopes}; it provides the prototype for the axisymmetric low-frequency result proved below.

\section{A Calder\'on-type inverse boundary value problem}\label{sec:calderon}
The Euclidean Laplacian in dimension $n$ is denoted by $\Delta_n$, and $\Delta_\Gamma$ denotes the surface Laplacian. Unmarked angle brackets are sesquilinear, linear in the first argument and conjugate-linear in the second; the subscript $\mathrm{bil}$ denotes a complex bilinear pairing without conjugation. Reciprocity identities use bilinear products, while energy and variational forms are sesquilinear.

\subsection{Magnetic responses as restricted electromagnetic boundary data}\label{subsec:restricted}
Fix $s\neq0$ with $\Re s\geq0$. For an incident magnetic-potential trace $f$, let $\widetilde f=\mathcal{Q}_\sigma(s)f$ be the induced trace and $\tau=f+\widetilde f$ the total trace. From \eqref{eq:gauss}, the boundary magnetic field is
\begin{equation}\label{eq:magtraces}
 \bfH_T=-\nabla_\Gamma\tau,\qquad
 RH_{\nu,\ell m}=-\ell f_{\ell m}+(\ell+1)\widetilde f_{\ell m}.
\end{equation}
Thus the magnetic input prescribed from the interior perspective is necessarily a surface-gradient tangential field. Define the corresponding magnetic boundary operator by
\begin{equation}
 \mathcal{T}_\sigma(s)\tau=-H_\nu,
\end{equation}
where $\bfH$ solves \eqref{eq:maxwell} with $\bfH_T=-\nabla_\Gamma\tau$. Let $\mathcal{D}_+$ and $\mathcal{D}_-$ be the Dirichlet--Neumann operators for the internally regular and externally decaying harmonic continuations, respectively. Their eigenvalues on degree $\ell$ are $\ell/R$ and $-(\ell+1)/R$, respectively; set $\mathcal{W}=\mathcal{D}_+-\mathcal{D}_-$, whose eigenvalues are $(2\ell+1)/R$. From $\tau=(I+\mathcal{Q}_\sigma)f$ and $\mathcal{T}_\sigma\tau=\mathcal{D}_+f+\mathcal{D}_-\mathcal{Q}_\sigma f$ one obtains
\begin{equation}\label{eq:QvsT}
 \mathcal{T}_\sigma=(\mathcal{D}_++\mathcal{D}_-\mathcal{Q}_\sigma)(I+\mathcal{Q}_\sigma)^{-1},\qquad
 \mathcal{Q}_\sigma=(\mathcal{T}_\sigma-\mathcal{D}_-)^{-1}(\mathcal{D}_+-\mathcal{T}_\sigma).
\end{equation}
These are operator identities between the appropriate zero-mean trace spaces and their duals, not merely entrywise transformations of a finite matrix.

The distinction from complete Maxwell boundary measurements can be seen directly. On the sphere a tangential vector field decomposes into surface-gradient and rotated-gradient components. The prescribed $\bfH_T$ is a gradient, whereas the normal component of Faraday's law gives
\begin{equation}\label{eq:surfCurl}
 \curl_\Gamma\bfE_T=\nu\cdot\curl\bfE=-s\mu H_\nu
 =s\mu \mathcal{T}_\sigma(s)\tau.
\end{equation}
If $\bfE_T=\nabla_\Gamma u+\nu\times\nabla_\Gamma v$, the surface curl determines $\Delta_\Gamma v$ but contains no information about the gradient part $\nabla_\Gamma u$. Accordingly, the data are restricted in \emph{both} their realizable boundary inputs and their measured electromagnetic output, even when available on the whole surface. This is different from having full Maxwell Cauchy data only on part of the boundary.

For completeness we verify the interior boundary map and invertibility in \eqref{eq:QvsT}. Eliminating $\bfE$ yields
\begin{equation}\label{eq:magBVP}
 \curl(\sigma^{-1}\curl\bfH)+s\mu\bfH=0,\qquad
 \bfH_T=-\nabla_\Gamma\tau.
\end{equation}
After prescribing the tangential trace, the variational form on $H_0(\curl)$ is coercive for $s\neq0$, $\Re s\geq0$, after multiplication by $e^{-i\arg(s)/2}$. The solution is divergence-free and has a normal trace. To check the boundary-map sign, on an oriented local frame $(\mathbf t_1,\mathbf t_2,\nu)$ one has $\divg_\Gamma(\nu\times\bfE)=-\curl_\Gamma\bfE_T$. Integration by parts, the identity $\divg(\bfE\times\overline{\bfH})=\overline{\bfH}\cdot\curl\bfE-\bfE\cdot\curl\overline{\bfH}$, and \eqref{eq:surfCurl} give
\begin{equation}\label{eq:Tpositive}
 \sesinner{\mathcal{T}_\sigma\tau}{\tau}
 =\norm{\bfH}_{L^2(\B)}^2
 +(s\mu)^{-1}\int_{\B}\sigma^{-1}|\curl\bfH|^2\dd x.
\end{equation}
The nonnegative real part of this form and the coercivity of $-\mathcal{D}_-$ imply that $\mathcal{T}_\sigma-\mathcal{D}_-$ is invertible from the zero-mean $H^{1/2}$ trace space to its dual. Since $\mathcal{W}$ is an isomorphism between the same spaces, $I+\mathcal{Q}_\sigma$ is invertible as well. Green formulas justify the weak-trace version for bounded $\sigma$.

The classical scalar Calder\'on map sends boundary potential to normal current for $\divg(\sigma\nabla u)=0$ \cite{SylvesterUhlmann1987}. Standard Maxwell inverse boundary value results use relations involving unrestricted tangential electric and magnetic traces \cite{Ola1993,Kenig2011}. Equation \eqref{eq:QvsT} instead identifies the particular \emph{restricted magnetic boundary operator} available to global induction; it does not equate $\mathcal{Q}_\sigma$ with the full Maxwell Cauchy-data map.

\subsection{An Alessandrini-type identity for the magnetic response}\label{subsec:bilinear}
The following identity plays the role of Alessandrini's identity in the classical Calder\'on problem. It converts a difference of boundary response operators into a volume integral weighted by the conductivity difference. In particular, equality of the complete magnetic responses implies an orthogonality relation between the conductivity difference and products of interior normalized electric fields. Throughout this subsection, reciprocity products are \emph{bilinear}, without complex conjugation.
\begin{proposition}[Alessandrini-type identity for the magnetic response]\label{prop:identity}
Fix $s\neq0$ with $\Re s\geq0$. For any admissible conductivity $\sigma$ and incident trace $f$, let $\eg^f_{\sigma,s}=\bfE^f_{\sigma,s}/(-\mu s)$ denote the normalized electric field, and define $\mathcal{F}_\sigma(s)=\mathcal{W}\mathcal{Q}_\sigma(s)/(\mu s)$. Then, for $\sigma_1,\sigma_2$ satisfying \eqref{eq:condbounds} and all incident traces $f,g$,
\begin{equation}\label{eq:alessandrini}
 \inner{[\mathcal{F}_{\sigma_1}(s)-\mathcal{F}_{\sigma_2}(s)]f}{g}
 =\int_{\B}(\sigma_1-\sigma_2)\,
 \eg_{\sigma_1,s}^{f}\cdot\eg_{\sigma_2,s}^{g}\dd x.
\end{equation}
The same identity extends to $s=0$ with the zero-frequency limits defined below.
\end{proposition}
\begin{proof}
Let $(\bfE_1,\bfH_1)$ and $(\bfE_2,\bfH_2)$ be fields for $\sigma_1$ and $\sigma_2$ with incident traces $f,g$. Maxwell's equations give
\begin{align}
 \divg(\bfE_1\times\bfH_2)&=-s\mu\bfH_1\cdot\bfH_2-\sigma_2\bfE_1\cdot\bfE_2,\\
 \divg(\bfE_2\times\bfH_1)&=-s\mu\bfH_1\cdot\bfH_2-\sigma_1\bfE_1\cdot\bfE_2.
\end{align}
Subtracting and applying the divergence theorem yields
\begin{equation}\label{eq:crossFlux}
 \int_{\B}(\sigma_1-\sigma_2)\bfE_1\cdot\bfE_2\dd x
 =\int_\Gamma(\bfE_1\times\bfH_2-\bfE_2\times\bfH_1)\cdot\nu\dd S.
\end{equation}
Write $\widetilde f=\mathcal{Q}_{\sigma_1}f$ and $\widetilde g=\mathcal{Q}_{\sigma_2}g$ for the induced traces, and set $\tau_1=f+\widetilde f$, $\tau_2=g+\widetilde g$. Since $\bfH_{2,T}=-\nabla_\Gamma\tau_2$ and $\curl_\Gamma\bfE_{1,T}=-s\mu H_{1,\nu}$, surface integration by parts gives
\begin{equation}
 \int_\Gamma(\bfE_1\times\bfH_2)\cdot\nu\dd S
 =-\inner{\curl_\Gamma\bfE_{1,T}}{\tau_2}
 =s\mu\inner{H_{1,\nu}}{\tau_2}.
\end{equation}
Interchanging indices and using $H_{1,\nu}=-\mathcal{D}_+f-\mathcal{D}_-\widetilde f$ and $H_{2,\nu}=-\mathcal{D}_+g-\mathcal{D}_-\widetilde g$, symmetry of $\mathcal D_\pm$ under bilinear pairing cancels the incident--incident and induced--induced terms:
\begin{align}
 \inner{H_{1,\nu}}{\tau_2}-\inner{H_{2,\nu}}{\tau_1}
 &=\inner{\mathcal{W}\widetilde f}{g}-\inner{\mathcal{W}\widetilde g}{f}.
\end{align}
For $\sigma_1=\sigma_2$, \eqref{eq:crossFlux} yields the reciprocity relation
\begin{equation}\label{eq:reciprocity}
 \inner{\mathcal{W}\mathcal{Q}_\sigma f}{g}=\inner{\mathcal{W}\mathcal{Q}_\sigma g}{f}.
\end{equation}
Thus, in a real spherical-harmonic basis, the matrix representing $\mathcal{W}\mathcal{Q}_\sigma$ is complex symmetric. Apply reciprocity to the second term for distinct conductivities, divide \eqref{eq:crossFlux} by $(\mu s)^2$, and use the definition of $\mathcal{F}_\sigma$ to obtain \eqref{eq:alessandrini}. All integrations extend to bounded conductivities using the $H(\curl)$ Green formula and weak boundary traces. The extension to zero frequency follows from \cref{subsec:zero_cont}.
\end{proof}

\subsection{Admissible currents and continuation to zero frequency}\label{subsec:zero_cont}
The low-frequency limit is most transparent in a current formulation. For $s\neq0$, the normalized induced current is
$\bfP_{\sigma,s}^{f}=\sigma\eg_{\sigma,s}^{f}=\bfJ_{\sigma,s}^{f}/(-\mu s)$.
Since $\divg\bfJ=0$ and the exterior is insulating, this current is divergence-free in $\B$ and has zero normal trace on $\Gamma$. The natural space of admissible normalized currents is therefore
\begin{equation}\label{eq:admissiblecurrent}
 X=\left\{\bfP\in L^2(\B)^3:
 \int_{\B}\bfP\cdot\nabla\eta\dd x=0
 \quad\text{for every }\eta\in H^1(\B)\right\}.
\end{equation}
The defining condition is the weak formulation of $\divg\bfP=0$ in $\B$ together with $\bfP\cdot\nu=0$ on $\Gamma$. Thus $X$ encodes current conservation and the insulating boundary condition, rather than an arbitrary auxiliary constraint.

Let $u_f$ be the harmonic continuation of an incident magnetic-potential trace to $\B$. For a solid harmonic $u_\ell$ of degree $\ell$, define
\begin{equation}\label{eq:vectorpot}
 \bfA_f=\frac{\xx\times\nabla u_\ell}{\ell+1},\qquad f=u_\ell|_\Gamma,
\end{equation}
then extend by spherical-harmonic linearity. Since
\begin{equation}
 \curl(\xx\times\nabla u_\ell)
 =\xx\Delta_3 u_\ell-3\nabla u_\ell+(\nabla u_\ell\cdot\nabla)\xx
 -(\xx\cdot\nabla)\nabla u_\ell=-(\ell+1)\nabla u_\ell,
\end{equation}
we have $\curl\bfA_f=-\nabla u_f$; $\bfA_f$ is divergence-free and tangent to the sphere. Write $\mathcal{V}f=\bfA_f$, defining the incident-field vector-potential map $\mathcal{V}$.

The space $X$ is closed in $L^2(\B)^3$. For $\bfP\in X$, its extension by zero to $\mathbb{R}^3$ is divergence-free in the distributional sense. Let $\mathcal N$ be the componentwise Newtonian potential of this zero extension:
\begin{equation}
 (\mathcal{N}\bfP)(x)=\frac1{4\pi}\int_{\B}\frac{\bfP(y)}{|x-y|}\dd y.
\end{equation}
The operator $\mathcal N$ is bounded on $L^2(\B)^3$. The magnetic field generated by the induced current is represented through the curl of this Newtonian potential. Faraday's law consequently gives
\begin{equation}
 \eg_{\sigma,s}^f=\bfA_f-s\mu \mathcal{N}\bfP_{\sigma,s}^f+\nabla v.
\end{equation}
The gradient correction is orthogonal to every element of $X$ by its defining weak condition, so projection onto $X$ eliminates it. Let $\mathcal{P}_X$ denote the $L^2$ orthogonal projection onto $X$, and define
\begin{equation}
 \mathcal{R}_\sigma=\mathcal{P}_X\sigma^{-1}|_X,\qquad \mathcal{L}=\mu \mathcal{P}_X\mathcal{N}|_X.
\end{equation}
Projecting the current-field relation onto $X$ then gives
\begin{equation}\label{eq:currentresolvent}
 \bfP_{\sigma,s}^f=(\mathcal{R}_\sigma+s\mathcal{L})^{-1}\mathcal{V}f.
\end{equation}
The quadratic form of $\mathcal{R}_\sigma$ is the ohmic dissipation, while that of $\mathcal{L}$ is the magnetic energy of the induced field:
$\sesinner{\mathcal{R}_\sigma\bfP}{\bfP}=\int_{\B}\sigma^{-1}|\bfP|^2\dd x$ and
$\sesinner{\mathcal{L}\bfP}{\bfP}=\mu\int_{\R^3}|\nabla \mathcal{N}\bfP|^2\dd x$.
Both forms are nonnegative. Moreover
\begin{equation}
 \Re\sesinner{(\mathcal{R}_\sigma+s\mathcal{L})\bfP}{\bfP}
 \geq\left[(\sigma^*)^{-1}+\min\{0,\Re s\}\norm{\mathcal{L}}\right]\norm{\bfP}_2^2.
\end{equation}
Thus \eqref{eq:currentresolvent} is analytic in a neighborhood of the closed right half-plane and $\norm{(\mathcal{R}_\sigma+s\mathcal{L})^{-1}}\leq\sigma^*$ for $\Re s\geq0$. This resolvent representation provides an $L^2$-analytic continuation of the normalized current, and hence of the normalized response $\mathcal F_\sigma(s)$, through $s=0$. Although $\mathcal Q_\sigma(0)=0$, the quotient $\mathcal W\mathcal Q_\sigma(s)/s$ therefore has a finite limit, which identifies the zero-frequency slope.

For a test trace $g$, integration of Amp\`ere's law against $\bfA_g$, using $\curl\bfA_g=-\nabla u_g$, leads to
\begin{equation}\label{eq:currentdata}
 \inner{\mathcal{F}_\sigma(s)f}{g}=\frac{\inner{\mathcal{W}\mathcal{Q}_\sigma(s)f}{g}}{\mu s}
 =\int_{\B}\bfP_{\sigma,s}^f\cdot\bfA_g\dd x.
\end{equation}
Indeed, for a test mode of degree $\ell$, the boundary contribution is
$-\frac{R}{\ell+1}\int_\Gamma\nabla_\Gamma\widetilde f\cdot\nabla_\Gamma g\dd S
=-\inner{\widetilde f}{\mathcal{D}_+g}$, and the volume contribution is $\inner{\mathcal{D}_-\widetilde f}{g}$. The identity follows by linearity and weak trace approximation. The right-hand side of \eqref{eq:currentdata} extends the matrix elements of $\mathcal{F}_\sigma$ analytically through $s=0$. Hence $\mathcal{Q}_\sigma(0)=0$, $\mathcal{W}\partial_s\mathcal{Q}_\sigma(0)=\mu \mathcal{F}_\sigma(0)$, and the normalized electric fields converge in $L^2$. Formula \eqref{eq:alessandrini} therefore extends to $s=0$. The analyticity underlying the Laplace and physical-frequency descriptions in \eqref{eq:physical-laplace} thus follows directly from the current resolvent and the current--response identity.

For later use, we record the precise equivalence between physical-frequency and analytic response data.
\begin{lemma}[Physical-frequency and analytic response data]\label{lem:frequency-equivalence}
Let $\sigma_1,\sigma_2$ satisfy \eqref{eq:condbounds}. Choose $\delta>0$ with $\delta\sigma^*\norm{\mathcal L}<1$ and set
$\mathcal U_\delta=\{s\in\mathbb C:\Re s>-\delta\}$.
For fixed input--output indices $(\ell,m)$ and $(k,n)$, put
$q_j(s)=Q_{\sigma_j;kn,\ell m}(s)$, $j=1,2$.
For any nonempty open interval $I\subset\mathbb R$,
\begin{equation}\label{eq:frequency-data-equivalence}
 \begin{aligned}
 q_1^{\mathrm{phys}}(\omega)=q_2^{\mathrm{phys}}(\omega)
 \quad\text{for every }\omega\in I
 \quad\Longleftrightarrow\quad
 q_1(s)=q_2(s)
 \quad\text{for every }s\in\mathcal U_\delta.
 \end{aligned}
\end{equation}
\end{lemma}
\begin{proof}
The coercivity bound for $\mathcal R_{\sigma_j}+s\mathcal L$ and \eqref{eq:currentdata} show that both entries are analytic in the common connected domain $\mathcal U_\delta$. Equality of their physical responses means that $q_1-q_2$ vanishes on $-iI$. This set has an accumulation point inside $\mathcal U_\delta$, so the identity theorem gives equality throughout that domain. The converse follows by restricting to $s=-i\omega$.
\end{proof}
Applying the lemma entrywise gives the same equivalence for any prescribed collection of spatial index pairs, including the complete matrix or the complete zonal block. Consequently, injectivity of the conductivity-to-data map is equivalent for exact physical-frequency curves on $I$ and for the corresponding analytic response family, with the spatial index set held fixed. Establishing either injectivity property still requires a conductivity-recovery argument; analyticity alone does not resolve the general three-dimensional problem.

\subsection{The three-dimensional uniqueness question and restricted product density}\label{subsec:threeD}
Let $s\neq0$ and suppose two conductivities have the same complete magnetic response. The Alessandrini-type identity then implies
\begin{equation}\label{eq:products}
 \int_{\B}(\sigma_1-\sigma_2)h\dd x=0,
 \qquad
 h\in\Span\left\{\eg_{\sigma_1,s}^{f}\cdot\eg_{\sigma_2,s}^{g}:\ f,g\text{ incident magnetic-potential traces}\right\}.
\end{equation}
Each product belongs to $L^1(\B)$. Density of this span in $L^1(\B)$ would be a sufficient condition for uniqueness among bounded conductivities. More generally, it is enough that the products separate all nonzero conductivity differences in the model class. Full $L^1$-density is therefore a convenient sufficient condition, but need not be a necessary characterization of uniqueness. Equation \eqref{eq:products} isolates the central issue as a \emph{density problem of the restricted solution-product form}.

This formulation is closely related to the standard strategy in Calder\'on-type inverse problems. In the scalar Calder\'on problem, after a Liouville transform, one constructs complex geometrical optics (CGO) solutions of the associated Schr\"odinger equation with leading behaviour $e^{x\cdot\zeta}(1+\psi_\zeta)$. Choosing two complex phases whose sum is a prescribed Fourier frequency makes the leading product carry the corresponding Fourier mode. Estimates for the remainder terms then allow recovery of the Fourier transform of the unknown coefficient. Related vector-valued CGO constructions are central to uniqueness proofs for the full Maxwell inverse boundary value problem \cite{SylvesterUhlmann1987,Ola1993,Kenig2011}.

Transferring this argument to the present setting requires control of the admissible boundary traces. Even when a CGO solution of the local Maxwell system is available, its boundary trace need not lie in the range of the admissible magnetic control operator, whose inputs are restricted to potential-field data. One must therefore establish a realization or approximation result for suitable CGO solutions, or directly construct admissible solution products with the required separating property. For comparison, realizable normalized electric fields satisfy
\begin{equation}
 \curl\curl\eg+s\mu\sigma\eg=0,
 \qquad\curl\eg=\bfH,
 \qquad\nu\cdot\sigma\eg=0,
\end{equation}
and their magnetic boundary control obeys
\begin{equation}
 (\curl\eg)\times\nu=\nu\times\nabla_\Gamma\tau.
\end{equation}
If $\mathbf z_T$ is the adjoint trace paired with this control, then
\begin{equation}\label{eq:RungeRestrict}
 \int_\Gamma\mathbf z_T\cdot(\nu\times\nabla_\Gamma\tau)\dd S
 =-\inner{\curl_\Gamma\mathbf z_T}{\tau}.
\end{equation}
Vanishing for every $\tau$ yields only $\curl_\Gamma\mathbf z_T=0$, which also holds for nonzero surface-gradient traces. It does not provide all the vanishing Cauchy data used in standard unrestricted vector Runge arguments. Even for constant conductivity, the accessible electric fields are toroidal, whereas general local Maxwell solutions may contain radial components. Crucially, density of \emph{scalar products} is weaker than density of all vector solutions. The adjoint boundary calculation above identifies an obstruction to applying the usual unrestricted Runge argument, but it is not evidence of non-uniqueness. Establishing a suitable restricted CGO or product-completeness result remains open, as does the general three-dimensional uniqueness or non-uniqueness of the complete magnetic $Q$-response.

\section{Axisymmetric uniqueness from the \texorpdfstring{$Q$}{Q}-response}\label{sec:axisym}
The general problem of \cref{subsec:threeD} is open, but the restricted data can be shown to be sufficient in several axisymmetric classes. The two main arguments are different. The zero-frequency slopes of the zonal response yield a complete family of \emph{linear moments} for arbitrary axisymmetric conductivity satisfying \eqref{eq:condbounds}. At a fixed nonzero frequency, a five-dimensional scalar reduction and a known shell structure yield \emph{boundary recovery and layer stripping by unique continuation and density of reachable traces}. We first establish the low-frequency result and its consequences, then the single-frequency result.

\subsection{Low-frequency moments and uniqueness from zonal slopes}\label{subsec:slopes}
\paragraph{A radial prototype.}
The relevance of the zero-frequency slope is already visible in the classical radially symmetric model. If $\sigma=\sigma(r)$, the response matrix is diagonal, with degree-$\ell$ entry $q_\ell(s)$, where $s=-i\omega$. Its first derivative has the explicit representation
\begin{equation}\label{eq:radialslope}
 q_\ell'(0)=\frac{\mu\ell}{(\ell+1)(2\ell+1)R^{2\ell+1}}
 \int_0^R\sigma(r)r^{2\ell+2}\dd r,\qquad \ell\geq1.
\end{equation}
The derivation from the radial induction equation and its relation to the $C_\ell$
response are given in \cref{app:radial}. With the change of variables $t=r^2$,
these moments become
\[
 \int_0^R\sigma(r)r^{2\ell+2}\,\dd r
 =\frac12\int_0^{R^2}t^{\ell-1}\sigma(\sqrt t)t^{3/2}\,\dd t.
\]
Thus the complete sequence of slopes gives all polynomial moments of the finite
measure $\frac12\sigma(\sqrt t)t^{3/2}\,\dd t$ on $[0,R^2]$. Polynomial
density in $C([0,R^2])$ implies uniqueness of this measure and hence of
$\sigma(r)$ almost everywhere. This spatial-moment mechanism is distinct from
the classical inverse-spectral recovery of a radial profile from the full
frequency response at one fixed degree. It also provides a simple prototype
for the two-dimensional moment mechanism used below.

\begin{remark}[Sparse degrees in the radial moment problem]
The radial moment argument does not require every degree. Let $\mathcal L\subset\{1,2,\ldots\}$ be nonempty, put $\ell_1=\min\mathcal L$, and suppose
\[
 \sum_{\substack{\ell\in\mathcal L\\\ell>\ell_1}}
 \frac{1}{\ell-\ell_1}=\infty.
\]
With $x=(r/R)^2$, the M\"untz--Sz\'asz theorem implies that
$\operatorname{span}\{x^{\ell-\ell_1}:\ell\in\mathcal L\}$ is dense in
$C([0,1])$. Hence the radial zero-frequency moments in~\eqref{eq:radialslope} associated with
only the degrees in $\mathcal L$ still suffice to determine $\sigma(r) r^{2\ell_1+2}$ (and thus $\sigma(r)$) almost everywhere. This is a statement about the one-variable
radial moment problem. Whether an analogous sparse selection of zonal
input--output degree pairs remains complete for general axisymmetric
conductivities is a separate multivariate approximation problem and is
not asserted here; see also \cite{BorweinErdelyi1996}.
\end{remark}

 Under axisymmetry, off-diagonal couplings between zonal degrees turn the radial moments into a complete family of moments in the two meridional variables.
 The following theorem extends this moment mechanism to axisymmetric conductivities.
\begin{theorem}[Uniqueness from zonal slopes]\label{thm:slope}
Let $\sigma_1,\sigma_2$ satisfy \eqref{eq:condbounds} and be rotationally symmetric about the same known axis. If
\begin{equation}\label{eq:slopeassump}
 \partial_\omega Q_{\sigma_1;k0,\ell0}(0)
 =\partial_\omega Q_{\sigma_2;k0,\ell0}(0)
 \quad\text{for all }k,\ell\geq1,
\end{equation}
then $\sigma_1=\sigma_2$ almost everywhere in $\B$.
\end{theorem}
The derivative is taken for each fixed conductivity. Since $s=-i\omega$, $\partial_\omega Q(0)=-i\partial_s Q(0)$. No smoothness beyond essential boundedness is assumed, and the conductivity contrast need not be small.

Let $\eg_{\sigma,0}^f=\sigma^{-1}\bfP_{\sigma,0}^f$. By \eqref{eq:currentresolvent}, $\eg_{\sigma,0}^f-\bfA_f$ is a gradient and $\sigma\eg_{\sigma,0}^f\in X$. We may therefore replace $\bfA_g$ by $\eg_{\sigma,0}^g$ in \eqref{eq:currentdata}. Define
\begin{equation}\label{eq:Kform}
 \frac1\mu\inner{\mathcal{W}\partial_s \mathcal{Q}_\sigma(0)f}{g}=K_\sigma(f,g),\qquad
 K_\sigma(f,g)=\int_{\B}\sigma\eg_{\sigma,0}^f\cdot\eg_{\sigma,0}^g\dd x.
\end{equation}
For zonal solid harmonics set $U_\ell=r^\ell P_\ell(z/r)$ and $f_\ell=U_\ell|_\Gamma$, where $P_\ell$ is a Legendre polynomial. Equation \eqref{eq:vectorpot} gives
\begin{equation}\label{eq:zonalvectorpot}
 \bfA_\ell=-\frac{\rho}{\ell+1}r^{\ell-1}P_\ell'(z/r)\ephi.
\end{equation}
This is an azimuthal vector field tangent to circles of fixed $(\rho,z)$. If $\sigma=\sigma(\rho,z)$, then for every smooth test $\eta$,
\begin{align}
 \int_{\B}\sigma\bfA_\ell\cdot\nabla\eta\dd x
 &=\int_{\substack{\rho>0\\\rho^2+z^2<R^2}}\sigma(\rho,z)A_\ell^\varphi(\rho,z)
 \left[\int_0^{2\pi}\partial_\varphi\eta(\rho,\varphi,z)\dd\varphi\right]\dd\rho\dd z=0.
\end{align}
By density this holds for every $\eta\in H^1(\B)$, hence $\sigma\bfA_\ell\in X$. At zero frequency, \eqref{eq:currentresolvent} then implies
\begin{equation}\label{eq:independence}
 \bfP_{\sigma,0}^{f_\ell}=\sigma\bfA_\ell,
 \qquad \eg_{\sigma,0}^{f_\ell}=\bfA_\ell.
\end{equation}
The normalized electric field is therefore \emph{independent of conductivity} in the axisymmetric zonal sector. At first order the induced current closes along azimuthal circles, without a conductivity-dependent gradient correction.

For $j\geq0$ put
\begin{equation}\label{eq:bj}
 b_j(\rho,z)=r^jP'_{j+1}(z/r)=r^jC_j^{(3/2)}(z/r),
\end{equation}
where $C_j^{(3/2)}$ is the Gegenbauer polynomial \cite[Ch.~22]{AbramowitzStegun1964}. The function $b_j$ is a homogeneous polynomial and depends on $\rho$ through $\rho^2$. From \eqref{eq:Kform} and \eqref{eq:independence},
\begin{equation}\label{eq:axisymmoments}
 K_\sigma(f_\ell,f_k)=\frac{1}{(\ell+1)(k+1)}
 \int_{\B}\sigma\rho^2 b_{\ell-1}b_{k-1}\dd x.
\end{equation}
Thus the first-order response is \emph{exactly linear} in $\sigma$, without any expansion in the conductivity contrast. With the unit-sphere normalization of $Y_{\ell0}$ specified in Section~\ref{sec:model}, put $K_{\ell k}=K_\sigma(f_\ell,f_k)$. Equation~\eqref{eq:Kform} gives
\begin{equation}\label{eq:moment-Q-normalization}
 \partial_s Q_{\sigma;k0,\ell0}(0)
 =\frac{\mu}{4\pi R^{k+\ell+1}}
  \sqrt{\frac{2\ell+1}{2k+1}}\,K_{\ell k},
 \qquad \partial_\omega Q(0)=-i\partial_s Q(0).
\end{equation}
In particular, the dipole coefficient satisfies
\begin{equation}\label{eq:dipole-moment}
 \partial_s Q_{\sigma;10,10}(0)
 =\frac{\mu}{16\pi R^3}\int_{\B}\sigma\rho^2\dd x.
\end{equation}
For constant conductivity this reduces to $\mu\sigma R^2/30$, consistently with \eqref{eq:radialslope}.

\begin{lemma}[Span of products of zonal factors]\label{lem:span}
For every integer $N\geq0$,
\begin{equation}\label{eq:polynomials}
 \Span\{b_jb_{N-j}:0\leq j\leq\lfloor N/2\rfloor\}
 =\Span\{\rho^{2p}z^{N-2p}:0\leq p\leq\lfloor N/2\rfloor\}.
\end{equation}
\end{lemma}
\begin{proof}
Write $C_j=C_j^{(3/2)}$. These polynomials are orthogonal on $[-1,1]$ for the weight $1-\xi^2$. Degree and parity allow expansion of $C_jC_{N-j}$ in the family $\{C_{N-2p}\}$. When $p>j$, the degree of $C_jC_{N-2p}$ is smaller than $N-j$, so it is orthogonal to $C_{N-j}$ and the corresponding coefficient vanishes. When $p=j$, the leading coefficient of $C_jC_{N-2j}$ is nonzero, hence its projection onto $C_{N-j}$ is nonzero. The expansion matrix is triangular with nonzero diagonal and therefore invertible.

Substituting $\xi=z/r$ and multiplying by $r^N$, it follows that the $b_jb_{N-j}$ span the family $r^{2p}b_{N-2p}$. All members lie in the polynomial space on the right of \eqref{eq:polynomials}. On restricting to $r=1$, they are polynomials of distinct degrees and hence linearly independent. The dimensions are equal, proving the identity.
\end{proof}
\begin{proof}[Proof of \cref{thm:slope}]
Let $h=\sigma_1-\sigma_2$. Equality of the slope matrices, together with \eqref{eq:Kform}, \eqref{eq:axisymmoments}, and \cref{lem:span}, yields
\begin{equation}
 \int_{\B}h(\rho,z)\rho^2P(\rho^2,z)\dd x=0
 \qquad\text{for every polynomial }P.
\end{equation}
Set $t=\rho^2$. After integrating over the azimuthal angle, $\dd x=\pi\dd t\dd z$ on the compact meridional set
\begin{equation}
 D=\{(t,z):t\geq0,\ t+z^2\leq R^2\}.
\end{equation}
The finite signed measure with density $t\,h(\sqrt t,z)$ has zero integral against every polynomial on $D$. Polynomials are uniformly dense in $C(D)$, so $t\,h(\sqrt t,z)=0$ almost everywhere. Since the symmetry axis has zero volume, $h=0$ almost everywhere in $\B$.
\end{proof}

\begin{corollary}[Open physical-frequency interval]\label{cor:interval}
Under the assumptions of \cref{thm:slope}, the exact complete zonal response block on any nonempty open interval of physical frequencies uniquely determines $\sigma$ almost everywhere.
\end{corollary}
\begin{proof}
Apply Lemma~\ref{lem:frequency-equivalence} to every zonal entry. Equality on the physical-frequency interval gives equality of the analytic entries and hence of their zero-frequency slopes. The conclusion follows from \cref{thm:slope}.
\end{proof}

\subsection{Cross-degree responses and two-dimensional moment completeness}\label{subsec:modeinfo}
The off-diagonal entries of the zonal response matrix are essential for recovering the full meridional dependence of an axisymmetric conductivity. Equation \eqref{eq:axisymmoments} shows that each pair of degrees $(k,\ell)$ supplies a moment weighted by the product $b_{k-1}b_{\ell-1}$. Diagonal entries alone have even parity in $z$ and therefore cannot distinguish a model from its reflection across the equatorial plane. Cross-degree entries with $k+\ell$ odd carry the corresponding odd moments. From $b_0=1$ and $b_1=3z$, for example,
\begin{equation}
 K_\sigma(f_1,f_1)=\frac14\int_{\B}\sigma\rho^2\dd x,
 \qquad K_\sigma(f_1,f_2)=\frac12\int_{\B}\sigma\rho^2z\dd x.
\end{equation}
The first is invariant under equatorial reflection, whereas the second changes sign. More generally, put $K_{\ell k}=K_\sigma(f_\ell,f_k)$ and
$M_{pq}=\int_{\B}\sigma\rho^{2p+2}z^q\dd x$. Since $b_2=6z^2-\frac32\rho^2$,
\begin{equation}\label{eq:momentsexamples}
 M_{00}=4K_{11},\quad M_{01}=2K_{12},\quad M_{02}=K_{22},
 \quad M_{10}=4K_{22}-\frac{16}{3}K_{13}.
\end{equation}
More generally, equatorial reflection multiplies the $(k,\ell)$-entry of the zonal Q-response by $(-1)^{k+\ell}$. Thus a model and its mirror image have the same diagonal responses at all frequencies, while odd cross-degree responses may differ. Beyond this parity information, the span identity in \cref{lem:span} shows that all degree-pair products together generate every polynomial in $(\rho^2,z)$. This is the mechanism by which the complete zonal slope matrix determines the two-dimensional meridional moments, rather than only radial or equatorially symmetric information.

If $\sigma(r,\theta)$ is, for almost every $r$, a polynomial of degree at most $L$ in $\xi=\cos\theta$, orthogonality with the weight $1-\xi^2$ shows that the zero-frequency slope entries with $|k-\ell|>L$ vanish. The banded infinite slope matrix is consequently sufficient for uniqueness in this class, but the allowed radial coefficient functions remain infinite-dimensional. A finite collection of slope entries contains only finitely many linear moments and cannot determine arbitrary bounded radial profiles. At nonzero frequency, the fields depend on $\sigma$, and the response need not retain the same bandwidth.

The slope theorem also implies pairwise discrimination for sufficiently low nonzero frequencies. For each fixed pair $\sigma_1\neq\sigma_2$, some slope entry differs, and
\begin{equation}
 Q_{\sigma_1;k0,\ell0}(\omega)-Q_{\sigma_2;k0,\ell0}(\omega)
 =\omega\,\partial_\omega(\mathcal{Q}_{\sigma_1}-\mathcal{Q}_{\sigma_2})_{k0,\ell0}(0)+O(\omega^2).
\end{equation}
That entry is unequal at all sufficiently small nonzero $\omega$, although the required threshold can depend on the pair. This statement concerns exact data; the response difference itself tends to zero as $\omega\to0$.

\paragraph{Finite slopes for finite-shell models.}
Assume known interfaces $0=r_0<r_1<\cdots<r_J=R$. Within shell $j$, let $\sigma$ be independent of $r$ and polynomial of degree at most $L$ in $\xi=z/r$. Expand the angular polynomial in $C_d^{(3/2)}(\xi)$ with coefficients $c_{jd}$. The basis $r^{2p}b_{N-2p}$ used in the proof of \cref{lem:span} allows construction of the tests $\rho^2r^{d+2n}C_d^{(3/2)}(\xi)$ from slope matrix entries with $\ell+k=d+2n+2$. The angular weight is $(1-\xi^2)C_d^{(3/2)}(\xi)$. After known normalization factors are removed, the data for each $d$ are
\begin{equation}\label{eq:finiteMoments}
 \sum_{j=1}^J c_{jd}\int_{r_{j-1}}^{r_j}r^{d+2n+4}\dd r,
 \qquad 0\leq n\leq J-1.
\end{equation}
For fixed $d$, the determinant of this $J\times J$ moment matrix is, up to a nonzero sign and constant,
\begin{equation}\label{eq:vandermonde}
 \int_{\prod_{j=1}^J(r_{j-1},r_j)}
 \left(\prod_{j=1}^J t_j^{d+4}\right)
 \prod_{1\leq j<k\leq J}(t_k^2-t_j^2)
 \dd t_1\cdots\dd t_J.
\end{equation}
The integrand is positive on the ordered product of shell intervals. All coefficients $c_{jd}$ are therefore determined. In particular, the zero-frequency slope submatrix with $1\leq k,\ell\leq L+2J-1$ suffices for this finite-dimensional model class.

\subsection{Single-frequency uniqueness for known concentric-shell models}
\label{subsec:shell}

We now consider a different uniqueness mechanism. The radial interfaces are assumed
known, while the conductivity may vary continuously with latitude within each shell.
Unlike the preceding low-frequency argument, the data in this subsection are taken at
one prescribed nonzero physical frequency.

\begin{theorem}[Single-frequency uniqueness for shell models]
\label{thm:shell}
Let $0=r_0<r_1<\cdots<r_J=R$ be known, and let two axisymmetric
conductivities satisfying \eqref{eq:condbounds} have the form
\begin{equation}\label{eq:shellclass}
 \sigma_i(r,\theta)=p_{i,j}(\cos\theta),\qquad
 r_{j-1}<r<r_j,\qquad
 p_{i,j}\in C([-1,1]),\qquad i=1,2.
\end{equation}
For any fixed $\omega_0\ne0$, equality of the complete zonal response blocks
\[
 \mathcal Q_{\sigma_1}^{00}(\omega_0)
 =
 \mathcal Q_{\sigma_2}^{00}(\omega_0)
\]
implies $\sigma_1=\sigma_2$ almost everywhere in $\mathbb{B}_R^3$.
\end{theorem}

\paragraph{Proof strategy.}
The proof has three stages. First, the azimuthal Maxwell branch generated by zonal
magnetic excitations is represented by an SO(4)-invariant scalar equation in a
five-dimensional ball. We show that the complete zonal magnetic response is equivalent
to the Dirichlet--Neumann map of this scalar equation restricted to the SO(4)-invariant
subspace. Second, high-degree invariant spherical harmonics extract all Fourier--cosine
coefficients of the conductivity in the outermost shell. Third, once an outer shell has
been identified, equality of the outer DN maps and unique continuation imply that the
two global solutions coincide in the common shell. A Runge-type density argument
shows that the traces of these global solutions are dense on the inner interface,
which then yields equality of the inner DN maps. Iterating this procedure over the
known interfaces proves the theorem.

\paragraph{Five-dimensional scalarization.}
Fix $s=-i\omega_0$, so that $s\ne0$ and $q=s\mu\sigma$ is purely imaginary.
For a zonal incident magnetic-potential trace, the axisymmetric Maxwell system has
an independent azimuthal electric-field branch. We write the normalized electric
field as
\begin{equation}\label{eq:toroidalfield}
 \eg=\rho U(\rho,z)\ephi,\qquad
 \rho=(x_1^2+x_2^2)^{1/2},\quad z=x_3.
\end{equation}
Here ``independent'' means that, for an axisymmetric coefficient and an $m=0$
excitation, the azimuthal component is decoupled from the meridional components in
the Maxwell operator. With the normalization
$\curl\eg=\bfH$ and $\bfE=-s\mu\eg$, Maxwell's equations reduce to
\begin{equation}\label{eq:curlcurl-normalized}
 \curl\curl\eg+s\mu\sigma\eg=0.
\end{equation}
A direct calculation gives
\begin{equation}\label{eq:curl-toroidal}
 \curl\eg=-\rho\,\partial_zU\,\mathbf e_\rho
 +(2U+\rho\,\partial_\rho U)\mathbf e_z
\end{equation}
and
\begin{equation}\label{eq:toroidal-scalar}
 -\partial_\rho^2U-\frac3\rho\partial_\rho U-\partial_z^2U
 +s\mu\sigma(\rho,z)U=0.
\end{equation}

To recognize the operator in \eqref{eq:toroidal-scalar}, introduce
$\mathbf y\in\mathbb{R}^4$ and set $\rho=|\mathbf y|$. We regard
\[
 U_5(\mathbf y,z)=U(|\mathbf y|,z).
\]
Thus $U_5$ is radial in the four variables $\mathbf y$, or equivalently invariant
under every rotation in $\mathrm{SO}(4)$. For such functions,
\[
 \Delta_4 U
 =\partial_\rho^2U+\frac3\rho\partial_\rho U.
\]
Consequently, \eqref{eq:toroidal-scalar} becomes the scalar equation
\begin{equation}\label{eq:5d}
 (-\Delta_5+q)U=0\quad\text{in }\mathbb{B}_R^5,\qquad
 q(\mathbf y,z)=s\mu\sigma(|\mathbf y|,z),
 \qquad U(\mathbf y,z)=U(|\mathbf y|,z).
\end{equation}
The five-dimensional problem is therefore an exact representation of the
axisymmetric toroidal Maxwell branch; it is not a change in the physical dimension
of the conductor. In particular, the original conductivity remains a function on
the three-dimensional ball.

The correspondence is also valid at the level of energy spaces. For smooth
SO(4)-invariant $U$, cylindrical integration gives
\begin{align}
 \|\rho U\ephi\|_{L^2(\mathbb{B}_R^3)}^2
 &=\frac1\pi\|U\|_{L^2(\mathbb{B}_R^5)}^2,\label{eq:5denergy1}\\
 \|\curl(\rho U\ephi)\|_{L^2(\mathbb{B}_R^3)}^2
 &=\frac1\pi\left(
 \|\nabla U\|_{L^2(\mathbb{B}_R^5)}^2+
 \frac2R\|U\|_{L^2(\partial \mathbb{B}_R^5)}^2\right).
 \label{eq:5denergy2}
\end{align}
These identities extend by density to the natural energy spaces. Their polarized
versions will be used to identify the weak Maxwell formulation with the five-dimensional
variational problem; details are given in Appendix~\ref{app:5d-energy}.

\paragraph{Boundary data and the invariant scalar DN map.}
Put
\[
 r=(\rho^2+z^2)^{1/2},\qquad
 \xi=\cos\theta=z/r,\qquad
 u=U|_{\partial \mathbb{B}_R^5}.
\]
From \eqref{eq:curl-toroidal}, the spherical components of the magnetic field are
\begin{equation}\label{eq:magnetic-boundary-components}
 H_\theta=-\sin\theta(2U+r\partial_rU),\qquad
 H_\nu|_{r=R}
 =-\partial_\xi\big[(1-\xi^2)u\big].
\end{equation}
The incident and total magnetic potentials determine the tangential magnetic trace
\[
 \bfH_T=-\nabla_\Gamma\tau.
\]
Matching the $\theta$-component in this identity gives the Robin condition
\begin{equation}\label{eq:5d-robin}
 \left(\partial_r+\frac2R\right)U
 =-\frac1{R^2}\partial_\xi\tau
 \quad\text{on }\partial \mathbb{B}_R^5.
\end{equation}
The normal magnetic output is
\begin{equation}\label{eq:5d-normal-output}
 H_\nu=-\partial_\xi\big[(1-\xi^2)u\big].
\end{equation}

We now formulate these relations as maps between boundary spaces. Let
$H_{\mathrm{zon},0}^t(\partial \mathbb{B}_R^3)$ denote the zero-mean zonal Sobolev
space on the three-dimensional sphere, and let
$H_{\mathrm{SO}(4)}^t(\partial \mathbb{B}_R^5)$ denote the subspace of
$H^t(\partial \mathbb{B}_R^5)$ invariant under rotations of the $\mathbf y$ variables.
Define, initially on finite harmonic expansions,
\begin{equation}\label{eq:DRMR}
 \mathcal D_R\tau=-R^{-2}\partial_\xi\tau,\qquad
 \mathcal M_Ru=-\partial_\xi\big[(1-\xi^2)u\big].
\end{equation}
The Legendre--Gegenbauer correspondence implies that these maps extend to
bounded isomorphisms
\begin{align}
 \mathcal D_R &:H_{\mathrm{zon},0}^{1/2}(\partial \mathbb{B}_R^3)
 \longrightarrow H_{\mathrm{SO}(4)}^{-1/2}(\partial \mathbb{B}_R^5),
 \label{eq:DRspaces}\\
 \mathcal M_R &:H_{\mathrm{SO}(4)}^{1/2}(\partial \mathbb{B}_R^5)
 \longrightarrow H_{\mathrm{zon},0}^{-1/2}(\partial \mathbb{B}_R^3).
 \label{eq:MRspaces}
\end{align}
The precise spectral calculation and Sobolev estimates are recorded in
Appendix~\ref{app:trace-isomorphisms}.

Let $\Lambda_q^{\mathrm{SO}(4)}$ be the Dirichlet--Neumann map for
\eqref{eq:5d}, restricted to the SO(4)-invariant subspace:
\begin{equation}\label{eq:invariant-DN}
 \Lambda_q^{\mathrm{SO}(4)}:
 H_{\mathrm{SO}(4)}^{1/2}(\partial \mathbb{B}_R^5)
 \longrightarrow H_{\mathrm{SO}(4)}^{-1/2}(\partial \mathbb{B}_R^5),
 \qquad
 \Lambda_q^{\mathrm{SO}(4)}u=\partial_\nu U|_{\partial \mathbb{B}_R^5}.
\end{equation}
Since $q$ is purely imaginary, the Dirichlet problem and the Robin problem
\eqref{eq:5d-robin} are well posed. Indeed, the real part of the Robin form is
\[
 \|\nabla U\|_{L^2(\mathbb{B}_R^5)}^2+
 \frac2R\|U\|_{L^2(\partial \mathbb{B}_R^5)}^2,
\]
which is coercive on $H^1(\mathbb{B}_R^5).$ If $u$ is the boundary trace of the Robin
solution, then
\begin{equation}\label{eq:Robin-DN-relation}
 \left(\Lambda_q^{\mathrm{SO}(4)}+\frac2R I\right)u
 =\mathcal D_R\tau.
\end{equation}

The remaining point is to verify that the Robin solution indeed generates the
three-dimensional Maxwell field with the prescribed tangential magnetic trace.
This is a weak-trace statement, not merely a pointwise calculation. It follows by
testing the five-dimensional variational equation against SO(4)-invariant scalar
functions and lifting them to three-dimensional azimuthal vector fields. The
polarized energy identities \eqref{eq:5denergy1}--\eqref{eq:5denergy2} give the
volume correspondence, while the Robin boundary term gives precisely the pairing
with $-\nabla_\Gamma\tau$. The normal output is then \eqref{eq:5d-normal-output}.
A complete verification, including the reduction from arbitrary vector tests by
rotation averaging, is given in Appendix~\ref{app:weak-maxwell-traces}.

Let $\mathcal T_\sigma^{00}$ be the zonal restriction of the magnetic boundary
operator in \eqref{eq:QvsT}. Equations \eqref{eq:5d-normal-output},
\eqref{eq:DRMR}, and \eqref{eq:Robin-DN-relation} yield
\begin{equation}\label{eq:Robin-Dirichlet-map}
 \left(\Lambda_q^{\mathrm{SO}(4)}+\frac2R I\right)^{-1}
 =-\mathcal M_R^{-1}\mathcal T_\sigma^{00}\mathcal D_R^{-1}.
\end{equation}
Here the identity operator is understood as the natural embedding from
$H_{\mathrm{SO}(4)}^{1/2}(\partial \mathbb{B}_R^5)$ into
$H_{\mathrm{SO}(4)}^{-1/2}(\partial \mathbb{B}_R^5)$.
Together with \eqref{eq:QvsT}, this proves the following result.

\begin{lemma}[Magnetic response and invariant scalar boundary map]
\label{lem:map_equivalence}
Let $\sigma$ satisfy \eqref{eq:condbounds} and be axisymmetric, and let
$\omega_0\ne0$ be fixed. With $s=-i\omega_0$ and
$q=s\mu\sigma$ lifted to $\mathbb{B}_R^5$ as in \eqref{eq:5d}, the complete zonal
response operator $\mathcal Q_\sigma^{00}(\omega_0)$ and the invariant
Dirichlet--Neumann map $\Lambda_q^{\mathrm{SO}(4)}$ determine each other.
\end{lemma}

\begin{proof}
The isomorphisms \eqref{eq:DRspaces}--\eqref{eq:MRspaces} and
\eqref{eq:Robin-Dirichlet-map} recover the Robin--Dirichlet map from
$\mathcal T_\sigma^{00}$ and conversely. The relation \eqref{eq:QvsT}
then gives the equivalence with $\mathcal Q_\sigma^{00}$.
\end{proof}

In particular,
\begin{equation}\label{eq:QimpliesLambda}
 \mathcal Q_{\sigma_1}^{00}(\omega_0)
 =\mathcal Q_{\sigma_2}^{00}(\omega_0)
 \quad\Longrightarrow\quad
 \Lambda_{q_1}^{\mathrm{SO}(4)}=\Lambda_{q_2}^{\mathrm{SO}(4)}.
\end{equation}

\paragraph{Recovery of the outermost shell.}
We next show that the invariant DN map determines the angular conductivity in
an outer shell. The underlying mechanism is concentration of high-degree
harmonic functions near the boundary. Pairing degrees whose difference remains
fixed retains the angular information in the limit.

\begin{lemma}[Determining the outermost shell]
\label{lem:outer}
Let $q$ in \eqref{eq:5d} be bounded and purely imaginary, and suppose that
\[
 q(r,\xi)=p^q(\xi)\qquad\text{for }a<r<R,
 \qquad p^q\in C([-1,1]).
\]
Set $x=R\widehat x$, $\widehat a=a/R$, and
$\widehat q(\widehat x)=R^2q(R\widehat x)$. The scaled DN map satisfies
$\widehat\Lambda_{\widehat q}=R\Lambda_q$, and the scaled outer-shell
profile is $\widehat p^q=R^2p^q$. Let $\Psi_n$ be the real,
$L^2$-normalized SO(4)-invariant spherical harmonics on $\mathbb{S}^4$,
with positive leading coefficient. Then, for every fixed integer $d\geq0$,
\begin{equation}\label{eq:highdegree}
 \lim_{n\to\infty}(2n+d+5)
 \left\langle
 (\widehat\Lambda_{\widehat q}^{\mathrm{SO}(4)}-\widehat\Lambda_0^{\mathrm{SO}(4)})\Psi_n,
 \Psi_{n+d}
 \right\rangle_{\mathrm{bil}}
 =
 \frac1\pi\int_0^\pi \widehat p^q(\cos\theta)\cos(d\theta)\dd\theta.
\end{equation}
Consequently, the invariant DN map determines $p^q$ and hence the
conductivity in the outer shell.
\end{lemma}

\begin{proof}
We work with the scaled variables and suppress hats. We divide the argument
into three steps.

\emph{Step 1: comparison with homogeneous harmonic functions.}
On the unit ball, set
\[
 v_n(r,\xi)=r^n\Psi_n(\xi).
\]
Then $v_n$ is harmonic, has boundary trace $\Psi_n$, and
\begin{equation}\label{eq:vn-norm}
 \|v_n\|_{L^2(\mathbb{B}_1^5)}=(2n+5)^{-1/2}.
\end{equation}
Let $u_n$ solve
\[
 (-\Delta_5+q)u_n=0\quad\text{in }\mathbb{B}_1^5,
 \qquad u_n|_{\mathbb{S}^4}=\Psi_n.
\]
The difference $u_n-v_n$ has zero Dirichlet trace and satisfies
\begin{equation}\label{eq:un-vn-resolvent}
 u_n-v_n=-(-\Delta_{5,D}+q)^{-1}(qv_n).
\end{equation}
The Dirichlet resolvent is compact on $L^2(\mathbb{B}_1^5)$: it maps $L^2$ boundedly
into $H_0^1$, and the latter embeds compactly into $L^2$. Since
$\{v_n/\|v_n\|_2\}$ is orthonormal, it converges weakly to zero. Multiplication
by $q\in L^\infty$ preserves weak convergence, and compactness in
\eqref{eq:un-vn-resolvent} gives
\begin{equation}\label{eq:un-vn-small}
 \|u_n-v_n\|_{L^2(\mathbb{B}_1^5)}=o(n^{-1/2}).
\end{equation}
The detailed compactness argument is given in Appendix~\ref{app:compact-resolvent}.

\emph{Step 2: reduction to a boundary integral.}
We use the bilinear Green identity (without complex conjugation). Since
$u_n$ solves the equation with potential $q$ and $v_{n+d}$ is harmonic,
\begin{equation}\label{eq:outerGreen}
 \left\langle
 (\Lambda_q^{\mathrm{SO}(4)}-\Lambda_0^{\mathrm{SO}(4)})\Psi_n,
 \Psi_{n+d}
 \right\rangle_{\mathrm{bil}}
 =
 \int_{\mathbb{B}_1^5}q\,u_n v_{n+d}\dd X.
\end{equation}
Replacing $u_n$ by $v_n$ incurs an error bounded by
\[
 \|q\|_\infty\|u_n-v_n\|_2\|v_{n+d}\|_2=o(n^{-1})
\]
for fixed $d$. Thus the right-hand side of \eqref{eq:outerGreen} equals
\[
 \int_{\mathbb{B}_1^5}q\,v_nv_{n+d}\dd X+o(n^{-1}).
\]
The contribution from $\mathbb{B}_{\widehat a}^5$ is bounded by
\[
 \|q\|_\infty
 \frac{\widehat a^{2n+d+5}}{2n+d+5},
\]
which is exponentially small after multiplication by $2n+d+5$.
In the scaled outer shell, $q=p^q(\xi)$, and radial integration gives
\begin{align}
 &(2n+d+5)\int_{\mathbb{B}_1^5\setminus \mathbb{B}_{\widehat a}^5}
 p^q(\xi)v_nv_{n+d}\dd X \notag\\
 &\qquad=(1-\widehat a^{2n+d+5})
 \int_{\mathbb{S}^4}p^q(\xi)\Psi_n(\xi)\Psi_{n+d}(\xi)\dd S.
 \label{eq:outer-shell-reduction}
\end{align}
Consequently, the limit in \eqref{eq:highdegree} reduces to the limit of the
spherical integral in \eqref{eq:outer-shell-reduction}.

\emph{Step 3: the Gegenbauer multiplication limit.}
The SO(4)-invariant spherical harmonics are normalized Gegenbauer polynomials
$\Psi_n(\xi)=c_nC_n^{(3/2)}(\xi)$. In this orthonormal basis, multiplication
by $\xi$ has the three-term recurrence
\begin{equation}\label{eq:gegenbauer-jacobi}
 \xi\Psi_n=a_n\Psi_{n+1}+a_{n-1}\Psi_{n-1},
 \qquad
 a_n=\sqrt{\frac{(n+1)(n+3)}{(2n+3)(2n+5)}}
 \longrightarrow\frac12.
\end{equation}
For a polynomial $p$, the matrix element
$\langle p(\xi)\Psi_n,\Psi_{n+d}\rangle$ is a finite sum of products of
the coefficients $a_{n+j}$, with bounded shifts $j$. As $n\to\infty$,
these coefficients converge to $1/2$. The limiting constant-coefficient
Jacobi operator has symbol
\[
 \frac{e^{i\theta}+e^{-i\theta}}2=\cos\theta.
\]
It follows that
\[
 \lim_{n\to\infty}
 \int_{\mathbb{S}^4}p(\xi)\Psi_n(\xi)\Psi_{n+d}(\xi)\dd S
 =
 \frac1\pi\int_0^\pi p(\cos\theta)\cos(d\theta)\dd\theta.
\]
Uniform polynomial approximation extends this identity to every
$p\in C([-1,1])$. The uniform error is bounded by the sup norm of the
polynomial approximation error, independently of $n$. Appendix~\ref{app:gegenbauer-limit}
contains the full Jacobi-matrix argument.

Combining the three steps proves \eqref{eq:highdegree}. The cosine
coefficients for all $d\geq0$ determine the continuous function
$p^q(\cos\theta)$ on $[0,\pi]$, and hence determine $p^q$ on $[-1,1]$.
Since the scaling factor and $s\mu$ are known and nonzero, the original
conductivity in the outer shell is determined.
\end{proof}

\paragraph{Removing a known shell.}
Once the outermost shell has been recovered, the two potentials coincide
there. We show that equality of the outer Dirichlet--Neumann maps then
implies equality of the Dirichlet--Neumann maps on the inner interface.
The argument has two ingredients: unique continuation in the common
shell, and density of the traces on the inner interface generated by
global solutions.

The key density property used below is the following. For a purely
imaginary, SO(4)-invariant potential $q\in L^\infty(\mathbb{B}_R^5)$, let $u_f$
solve the global Dirichlet problem with invariant boundary value
$f\in H_{\mathrm{SO}(4)}^{1/2}(\partial \mathbb{B}_R^5)$. Then
\begin{equation}\label{eq:reachable-trace-density}
 \overline{\left\{
 u_f|_{\partial \mathbb{B}_a^5}:
 f\in H_{\mathrm{SO}(4)}^{1/2}(\partial \mathbb{B}_R^5)
 \right\}}^{\,H^{1/2}(\partial \mathbb{B}_a^5)}
 =H_{\mathrm{SO}(4)}^{1/2}(\partial \mathbb{B}_a^5).
\end{equation}
The proof is by duality, using an adjoint interface-source problem and
unique continuation; the complete weak formulation is given in
Appendix~\ref{app:layer-stripping}.

\begin{lemma}[Removing a known shell]
\label{lem:strip}
Let $q_1,q_2\in L^\infty(\mathbb{B}_R^5)$ be purely imaginary and SO(4)-invariant,
and suppose that $q_1=q_2$ in
$\mathbb{B}_R^5\setminus\overline{\mathbb{B}_a^5}$ for some $0<a<R$. If their SO(4)-invariant
Dirichlet--Neumann maps on $\partial \mathbb{B}_R^5$ coincide, then their
corresponding invariant Dirichlet--Neumann maps on $\partial \mathbb{B}_a^5$
coincide.
\end{lemma}

\begin{proof}[Proof of \cref{lem:strip}]
Write $\Omega_{\mathrm{sh}}=\mathbb{B}_R^5\setminus\overline{\mathbb{B}_a^5}$ for the
common shell. Fix an arbitrary outer Dirichlet trace
$f\in H_{\mathrm{SO}(4)}^{1/2}(\partial \mathbb{B}_R^5)$, and let $u_i$ solve
\[
 (-\Delta_5+q_i)u_i=0\quad\text{in }\mathbb{B}_R^5,
 \qquad u_i|_{\partial \mathbb{B}_R^5}=f.
\]
The equality of the outer invariant DN maps gives equality of the outer
Neumann traces. Since the coefficients and the solutions are
SO(4)-invariant, this equality in the invariant boundary pairing is
equality of the full invariant Neumann traces. Consequently, for
$w=u_1-u_2$,
\[
 (-\Delta_5+q_{\mathrm{sh}})w=0\quad\text{in }\Omega_{\mathrm{sh}},
 \qquad
 w|_{\partial \mathbb{B}_R^5}=0,
 \quad
 \partial_\nu w|_{\partial \mathbb{B}_R^5}=0.
\]
Because both the Dirichlet and Neumann traces of $w$ vanish on the
outer sphere, its extension by zero across that sphere is a weak
solution of the equation with the common bounded potential extended
by zero. The unique-continuation theorem for Schr\"odinger operators with
bounded potentials (see, e.g., \cite{JerisonKenig1985}) then yields
\begin{equation}\label{eq:solutions-equal-shell}
 u_1=u_2\quad\text{in }\Omega_{\mathrm{sh}}.
\end{equation}
In particular, their traces on the inner interface agree. Denote the
common trace by
\[
 t_f=u_1|_{\partial \mathbb{B}_a^5}=u_2|_{\partial \mathbb{B}_a^5}.
\]
The weak transmission condition across $\partial \mathbb{B}_a^5$ identifies the
interior normal derivative, taken with the outward normal of $\mathbb{B}_a^5$,
with the derivative of the shell solution in the same geometric normal
direction. (The outward normal of the shell on its inner boundary is
$-\nu_a$, where $\nu_a$ is the outward normal of $\mathbb{B}_a^5$.) Since
the two shell solutions in \eqref{eq:solutions-equal-shell} coincide,
their derivatives in the direction $\nu_a$ coincide. Thus
\begin{equation}\label{eq:inner-DN-on-reachable}
 (\Lambda_1-\Lambda_2)t_f=0
 \qquad\text{for every }f\in
 H_{\mathrm{SO}(4)}^{1/2}(\partial \mathbb{B}_R^5),
\end{equation}
where $\Lambda_i$ denotes the invariant DN map of the inner ball with
potential $q_i|_{\mathbb{B}_a^5}$.

By the density property \eqref{eq:reachable-trace-density}, the set of traces $t_f$ is dense in
$H_{\mathrm{SO}(4)}^{1/2}(\partial \mathbb{B}_a^5)$. Since both inner DN maps are
bounded from $H_{\mathrm{SO}(4)}^{1/2}$ to $H_{\mathrm{SO}(4)}^{-1/2}$,
\eqref{eq:inner-DN-on-reachable} extends by continuity to every
admissible inner Dirichlet trace. Hence
\[
 \Lambda_1=\Lambda_2,
\]
which proves the lemma.
\end{proof}

\begin{proof}[Proof of \cref{thm:shell}]
By \eqref{eq:QimpliesLambda}, equality of the complete zonal Q-response
blocks gives equality of the invariant scalar DN maps on $\partial \mathbb{B}_R^5$.
Lemma~\ref{lem:outer} determines the conductivity in the outermost
shell. Lemma~\ref{lem:strip} then gives equality of the invariant scalar
DN maps on the next inner sphere. Repeating the outer-shell recovery and
shell-removal steps finitely many times over the known interfaces
determines every $p_{i,j}$. Since $q_i=s\mu\sigma_i$ and $s\mu\ne0$,
we conclude that $\sigma_1=\sigma_2$ almost everywhere in $\mathbb{B}_R^3$.
\end{proof}

\subsection{Finite zonal data for finite-dimensional shell models}\label{subsec:finite}
Returning to the low-frequency moment mechanism, we now consider a finite-dimensional shell prior. Let $0=r_0<r_1<\cdots<r_J=R$ be the known interfaces and write $\Omega_j=\{x:r_{j-1}<|x|<r_j\}$. Define the finite-dimensional difference space
\begin{equation}\label{eq:finiteModelSpace}
 V_{J,L}=\left\{h:\ h|_{\Omega_j}=p_j(\cos\theta),\quad
 p_j\in\mathbb P_L,\quad j=1,\ldots,J\right\},
\end{equation}
where $\mathbb P_L$ denotes the real polynomials of degree at most $L$. Differences of two conductivities in the corresponding shell model class belong to $V_{J,L}$.

The zero-frequency slope data define a finite moment map on this space. The moment calculation in \cref{subsec:modeinfo} shows that this map is injective when the maximum zonal degree is chosen as $N=L+2J-1$. Finite dimensionality then gives a positive lower bound for the map on the unit sphere of $V_{J,L}$. The resolvent estimates from \cref{subsec:zero_cont} show that the normalized response at a small nonzero frequency is a uniform $O(|s|)$ perturbation of this moment map. This yields the finite-frequency uniqueness and stability estimate below.

\begin{corollary}[Finite-mode uniqueness and stability]\label{cor:finite}
Fix $\mu$, the known interfaces $0=r_0<\cdots<r_J=R$, an integer $L\geq0$, and the bounds $\sigma_*,\sigma^*$. Let $\sigma$ satisfy \eqref{eq:condbounds} and, in each shell $\Omega_j$, have the form $p_j(\cos\theta)$ with $p_j\in\mathbb P_L$. Put $N=L+2J-1$ and
\begin{equation}
 Q_{\sigma,N}^{00}(\omega)=\bigl[Q_{\sigma;k0,\ell0}(\omega)\bigr]_{1\leq k,\ell\leq N}.
\end{equation}
There exist $\omega_*>0$ and $C_*>0$ depending only on the stated prior parameters such that for every $0<|\omega|<\omega_*$ and any two models in the class,
\begin{equation}\label{eq:stability}
 \norm{\sigma_1-\sigma_2}_{L^\infty(\B)}
 \leq\frac{C_*}{|\omega|}
 \norm{Q_{\sigma_1,N}^{00}(\omega)-Q_{\sigma_2,N}^{00}(\omega)}_F.
\end{equation}
Thus a finite zonal submatrix at a single such frequency uniquely determines the model.
\end{corollary}
\begin{proof}
Use real $L^2(\partial\B)$-normalized zonal inputs $\widehat f_\ell(R\widehat x)=R^{-1}Y_{\ell0}(\widehat x)$, and write $\widehat{\bfA}_\ell=\mathcal{V}\widehat f_\ell$. Multiplying every input and output basis function by the same factor $R^{-1}$ leaves the response matrix entries unchanged. Let $F_{\sigma,N}$ be the first $N\times N$ block in \eqref{eq:currentdata}. The finite-slope moment map is
$\mathcal{S}_N:V_{J,L}\to\mathbb{R}^{N\times N}$, defined by the formula below. At $s=0$,
\begin{equation}\label{eq:SN}
 F_{\sigma_1,N}(0)-F_{\sigma_2,N}(0)=\mathcal{S}_N(\sigma_1-\sigma_2),
 \qquad (\mathcal{S}_Nh)_{k\ell}=\int_{\B}h\widehat{\bfA}_\ell\cdot\widehat{\bfA}_k\dd x.
\end{equation}
The determinant argument in \cref{subsec:modeinfo} proves that $\mathcal{S}_N$ is injective on $V_{J,L}$. By finite dimensionality,
\begin{equation}\label{eq:kappa}
 \kappa=\min_{\substack{h\in V_{J,L}\\\norm h_\infty=1}}\norm{\mathcal{S}_Nh}_F>0.
\end{equation}

Write $\mathcal{G}_\sigma(s)=(\mathcal{R}_\sigma+s\mathcal{L})^{-1}$. For $\Re s\geq0$, the positivity bounds give $\norm{\mathcal{G}_\sigma(s)}\leq\sigma^*$. Since $\eg_{\sigma,0}^{\widehat f_\ell}=\widehat{\bfA}_\ell$, the resolvent identity yields
\begin{equation}\label{eq:resolventbounds}
 \norm{\eg_{\sigma,s}^{\widehat f_\ell}}_2
 \leq\frac{\sigma^*}{\sigma_*}\norm{\widehat{\bfA}_\ell}_2,
 \qquad
 \norm{\eg_{\sigma,s}^{\widehat f_\ell}-\widehat{\bfA}_\ell}_2
 \leq |s|\frac{(\sigma^*)^2}{\sigma_*}\norm{\mathcal{L}}\norm{\widehat{\bfA}_\ell}_2.
\end{equation}
These bounds are uniform over the model class.

Let $h=\sigma_1-\sigma_2$ and $\Delta F_N(s)=F_{\sigma_1,N}(s)-F_{\sigma_2,N}(s)$. By \eqref{eq:alessandrini}, the $(k,\ell)$ entry of $\Delta F_N(s)-\mathcal{S}_Nh$ is
\begin{equation}
 \int_{\B}h\left[
 (\eg_{\sigma_1,s}^{\widehat f_\ell}-\widehat{\bfA}_\ell)\cdot\eg_{\sigma_2,s}^{\widehat f_k}
 +\widehat{\bfA}_\ell\cdot(\eg_{\sigma_2,s}^{\widehat f_k}-\widehat{\bfA}_k)\right]\dd x.
\end{equation}
Cauchy--Schwarz and \eqref{eq:resolventbounds}, followed by the Frobenius norm, give a uniform perturbation estimate
\begin{equation}\label{eq:DeltaF}
 \norm{\Delta F_N(s)-\mathcal{S}_Nh}_F\leq C_N|s|\norm h_\infty,
 \qquad
 C_N=\frac{(\sigma^*)^2}{\sigma_*}\norm{\mathcal{L}}\left(1+\frac{\sigma^*}{\sigma_*}\right)
 \sum_{\ell=1}^N\norm{\widehat{\bfA}_\ell}_2^2.
\end{equation}
For $s=-i\omega$ with $0<|\omega|<\omega_*=\kappa/(2C_N)$, it follows that
$\norm{\Delta F_N(-i\omega)}_F\geq(\kappa/2)\norm h_\infty$. Equation \eqref{eq:currentdata} gives
\begin{equation}
 \Delta F_N(-i\omega)
 =\frac{W_N}{-i\mu\omega}
 \left[Q_{\sigma_1,N}^{00}(\omega)-Q_{\sigma_2,N}^{00}(\omega)\right],
 \quad W_N=\operatorname{diag}_{1\leq k\leq N}\frac{2k+1}{R}.
\end{equation}
Taking $C_*=2\norm{W_N}/(\mu\kappa)$ proves \eqref{eq:stability}.
\end{proof}
The estimate concerns exact $Q$ data at the chosen frequency. The factor $1/|\omega|$ increases error amplification as the frequency decreases, and the threshold $\omega_*$ may depend unfavorably on the finite-moment conditioning $\kappa$ as $J$ or $L$ increases. Quantitative bounds on these constants are not claimed here.

\section{Discussion}\label{sec:discussion}
The low-frequency analysis provides an explicit interpretation of the conductivity information encoded in magnetic $Q$-responses. For an axisymmetric conductivity under zonal excitation, the leading induced current follows azimuthal circles and is proportional to $\sigma\bfA_\ell$, where $\bfA_\ell$ is determined by the incident mode. This structure makes the zero-frequency slope a linear functional of conductivity, even when the conductivity contrast is large. The moment representation therefore complements numerical inversion by identifying the spatial quantities measured by individual response slopes.

With the normalization used in Section~\ref{subsec:slopes}, the associated bilinear form is
\[
 K_\sigma(f_\ell,f_k)
 =\frac{1}{(\ell+1)(k+1)}
 \int_{\B}\sigma\rho^2b_{\ell-1}b_{k-1}\,\dd x.
\]
Each pair of incident and induced degrees selects a polynomial weight in the meridional coordinates, with the relation to the measured response specified by \eqref{eq:moment-Q-normalization}. Even when the complete slope matrix is unavailable, individual slopes retain an explicit interpretation as conductivity-weighted spatial integrals, without requiring reconstruction of the entire conductivity distribution. Completeness of the full moment family gives the uniqueness theorem. At finite nonzero frequency, the internal fields depend on conductivity, and the response generally depends nonlinearly on it.

The off-diagonal coefficients are particularly informative about lateral structure. Cross-degree terms provide mixed moments, including moments odd in $z$, and therefore carry information about north--south asymmetry that is absent from the diagonal responses. For example, the degree-one and degree-two slope coupling measures a quantity proportional to $\int_{\B}\sigma\rho^2z\,\dd x$. A nonzero value detects equatorial asymmetry, while a zero value does not establish symmetry, because contributions from different regions may cancel. Each coefficient thus constrains a global weighted integral; it does not by itself locate a conductivity anomaly.

The same bilinear form has an energy interpretation: its diagonal entries determine the leading coefficients of low-frequency Joule dissipation, up to the known excitation normalizations, while its off-diagonal entries describe the cross terms for superposed real excitations. Under suitable integrability assumptions, the zero-frequency slope is also related to the first temporal moment of a causal response, providing a connection with the mantle-filter perspective of Backus \cite{Backus1983}. Since the induced-to-incident $Q$-response has zero static gain, this moment does not define a normalized delay time.

The spatial weights also reveal limitations of recovery. For conductivity perturbations confined to $r<a<R$, homogeneity of the vector potentials gives a moment bound with the attenuation factor $(a/R)^{k+\ell+3}$ relative to their full-ball $L^2$ norms. Deep contributions are therefore increasingly suppressed at high degree. The boundary concentration used in the single-frequency argument is consistent with this behavior, and the layer-stripping proof supplies no uniform stability estimate for unrestricted continuous shell profiles. In the finite-dimensional result, the factor $1/|\omega|$ in \eqref{eq:stability} exposes a further tradeoff: lower frequencies improve the leading-order approximation while increasing the sensitivity of recovery to response errors.

Observational access imposes additional restrictions. Only finitely many spatial modes are available in practice, and their estimation depends on source diversity, source separation, measurement noise, and continuation from the observation surface. Estimating a $Q$-response matrix and recovering conductivity from it involve distinct stability questions. Within the respective axisymmetric model classes considered here, the uniqueness results establish ideal information benchmarks; practical resolution also depends on the spatial weights, accessible modes, and error levels.

\section{Outlook}\label{sec:outlook}
\subsection{General three-dimensional uniqueness}\label{subsec:open3d}
A central question left open by this work is whether the complete magnetic $Q$-response at a prescribed nonzero frequency determines a general three-dimensional scalar conductivity satisfying \eqref{eq:condbounds} in the known ball. The magnetic potential-field excitations and observations provide restricted boundary information, so the standard full-data Maxwell uniqueness arguments do not apply directly. Extending the integral-identity approach requires a suitable completeness result for products of the accessible fields, or another mechanism that uses the restricted data.

An intermediate question is single-frequency uniqueness for general axisymmetric conductivities without the shell assumption. The scalar reduction remains available, but the present layer-stripping argument relies on radial constancy within known shells and does not directly extend to arbitrary radial variation.

\subsection{Beyond spherical geometry}\label{subsec:generaldomains}
A second direction is magnetic-only conductivity recovery in a known, non-spherical conducting domain. Spherical geometry provides both the harmonic representation of the $Q$-response and the explicit fields underlying the moment argument. In a general domain, one may instead formulate an exterior magnetic input--output operator and investigate whether it determines the conductivity. This requires specifying the admissible excitations and measurement traces, accounting for the topology of the exterior region, and relating observations on a separate sensor surface to data at the conductor boundary.

The main question is which features of the restricted boundary formulation survive this change of geometry and support a corresponding uniqueness argument. Recovery of an unknown conductor shape or conductivity support is a related, separate problem: identifying the conducting region alone does not determine the conductivity distribution within it.

\subsection{Multi-frequency information and inversion}\label{subsec:multifrequency}
A third direction concerns the effective use of responses measured at several temporal frequencies. Geomagnetic $Q$-responses are commonly estimated at discrete frequencies within a limited band and used jointly in conductivity inversion \cite{Puthe2014,Kuvshinov2021}. The present results primarily exploit spatial-mode diversity through zero-frequency slopes or responses at one frequency. Although analyticity also gives uniqueness from an exact complete zonal response on an open frequency interval, this conclusion supplies no stable recovery procedure for finitely sampled, noisy data with limited spatial modes.

The low-frequency expansion identifies a useful starting point. For a fixed axisymmetric conductivity, let $Q_N(s)$ denote a fixed finite zonal block in the Laplace parameterization of Section~\ref{subsec:complex-frequency}, and set $M_N(\sigma)=\partial_s Q_N(0)$. Its physical response satisfies
\begin{equation}\label{eq:multifrequency-expansion}
 Q_N^{\mathrm{phys}}(\omega)=Q_N(-i\omega)=-i\omega M_N(\sigma)+O(\omega^2),
 \qquad \omega\to0.
\end{equation}
Additional low-frequency samples may improve estimation of $M_N(\sigma)$, but they do not introduce new spatial weights at first order. Higher-order frequency dependence may carry conductivity information beyond that encoded in the leading moment matrix. The radial inverse-spectral benchmark shows that the complete, exact frequency response of a single spatial mode can determine radial structure; whether comparable spatial--frequency complementarity holds in broader model classes requires further analysis.

This motivates studying how frequency diversity changes identifiability and stability for a prescribed set of spatial modes. For finite-dimensional models, useful questions include whether several frequencies permit recovery with fewer spatial coefficients and how frequency selection affects conditioning. Such complementarity has limits: an axisymmetric conductivity and its equatorial reflection have identical diagonal zonal responses at every frequency, so frequency diversity alone cannot remove this ambiguity. A practical theory should characterize the information supplied jointly by spatial coupling and temporal frequency, and develop ways to use the observed band directly while balancing noise, low-frequency approximation error, and sensitivity to different parts of the conductor.

\printbibliography[title={References}]

@article{Bailey1970,
  author  = {Bailey, R. C.},
  title   = {Inversion of the geomagnetic induction problem},
  journal = {Proceedings of the Royal Society of London A},
  volume  = {315},
  pages   = {185--194},
  year    = {1970},
  doi     = {10.1098/rspa.1970.0036}
}

@article{Backus1983,
  author  = {Backus, G. E.},
  title   = {Application of mantle filter theory to the magnetic jerk of 1969},
  journal = {Geophysical Journal International},
  volume  = {74},
  pages   = {713--746},
  year    = {1983},
  doi     = {10.1111/j.1365-246X.1983.tb01901.x}
}

@article{BorweinErdelyi1996,
  author  = {Borwein, P. and Erd{\'e}lyi, T.},
  title   = {The full {M\"untz} theorem in {\(C[0, 1]\)} and {\(L_1[0, 1]\)}},
  journal = {Journal of the London Mathematical Society},
  volume  = {54},
  pages   = {102--110},
  year    = {1996},
  doi     = {10.1112/jlms/54.1.102}
}

@article{Baratchart2021,
  author  = {Baratchart, L. and Gerhards, C. and Kegeles, A.},
  title   = {Decomposition of {\(L^2\)}-vector fields on {Lipschitz} surfaces: characterization via null-spaces of the scalar potential},
  journal = {SIAM Journal on Mathematical Analysis},
  volume  = {53},
  pages   = {4096--4117},
  year    = {2021},
  doi     = {10.1137/20M1387754}
}

@book{ChaveJones2012,
  editor    = {Chave, A. D. and Jones, A. G.},
  title     = {The Magnetotelluric Method: Theory and Practice},
  publisher = {Cambridge University Press},
  year      = {2012},
  doi       = {10.1017/CBO9781139020138}
}

@article{EckhardtKostenko2016,
  author  = {Eckhardt, J. and Kostenko, A.},
  title   = {The inverse spectral problem for indefinite strings},
  journal = {Inventiones Mathematicae},
  volume  = {204},
  number  = {3},
  pages   = {939--977},
  year    = {2016},
  doi     = {10.1007/s00222-015-0629-1}
}

@article{Glassmeier2014,
  author  = {Glassmeier, K.-H. and Tsurutani, B. T.},
  title   = {Carl Friedrich {Gauss} -- General Theory of Terrestrial Magnetism -- a revised translation of the {German} text},
  journal = {History of Geo- and Space Sciences},
  volume  = {5},
  pages   = {11--62},
  year    = {2014},
  doi     = {10.5194/hgss-5-11-2014}
}

@article{HolthamOldenburg2010,
  author  = {Holtham, E. and Oldenburg, D. W.},
  title   = {Three-dimensional inversion of {ZTEM} data},
  journal = {Geophysical Journal International},
  volume  = {182},
  pages   = {168--182},
  year    = {2010},
  doi     = {10.1111/j.1365-246X.2010.04634.x}
}

@article{HooshyarRazavy1982,
  author  = {Hooshyar, M. A. and Razavy, M.},
  title   = {The inverse problem of geomagnetic induction at a fixed frequency},
  journal = {Geophysical Journal International},
  volume  = {71},
  pages   = {139--150},
  year    = {1982},
  doi     = {10.1111/j.1365-246X.1982.tb04988.x}
}

@article{Huang2026,
  author  = {Huang, X. and Gerhards, C. and Ren, Z.},
  title   = {Global versus regional internal-external potential field separation},
  journal = {arXiv preprint arXiv:2603.04886},
  year    = {2026},
  doi     = {10.48550/arXiv.2603.04886}
}

@article{JerisonKenig1985,
  author  = {Jerison, D. and Kenig, C. E.},
  title   = {Unique continuation and absence of positive eigenvalues for {Schr\"odinger} operators},
  journal = {Annals of Mathematics},
  volume  = {121},
  pages   = {463--488},
  year    = {1985},
  doi     = {10.2307/1971205}
}

@article{Kenig2011,
  author  = {Kenig, C. E. and Salo, M. and Uhlmann, G.},
  title   = {Inverse problems for the anisotropic {Maxwell} equations},
  journal = {Duke Mathematical Journal},
  volume  = {157},
  pages   = {369--419},
  year    = {2011},
  doi     = {10.1215/00127094-1272903}
}

@article{Kuvshinov2021,
  author  = {Kuvshinov, A. and Grayver, A. and T{\o}ffner-Clausen, L. and Olsen, N.},
  title   = {Probing 3-{D} electrical conductivity of the mantle using 6 years of {Swarm}, {CryoSat-2} and observatory magnetic data and exploiting matrix {\(Q\)}-responses approach},
  journal = {Earth, Planets and Space},
  volume  = {73},
  pages   = {67},
  year    = {2021},
  doi     = {10.1186/s40623-020-01341-9}
}

@article{Mittelholz2026,
  author  = {Mittelholz, A.},
  title   = {A review of electromagnetic induction on planetary bodies},
  journal = {Surveys in Geophysics},
  volume  = {47},
  pages   = {619--675},
  year    = {2026},
  doi     = {10.1007/s10712-025-09896-6}
}

@book{AbramowitzStegun1964,
  editor    = {Abramowitz, M. and Stegun, I. A.},
  title     = {Handbook of Mathematical Functions with Formulas, Graphs, and Mathematical Tables},
  series    = {Applied Mathematics Series},
  number    = {55},
  publisher = {National Bureau of Standards},
  location  = {Washington, DC},
  year      = {1964}
}

@article{Ola1993,
  author  = {Ola, P. and P{\"a}iv{\"a}rinta, L. and Somersalo, E.},
  title   = {An inverse boundary value problem in electrodynamics},
  journal = {Duke Mathematical Journal},
  volume  = {70},
  pages   = {617--653},
  year    = {1993},
  doi     = {10.1215/S0012-7094-93-07014-7}
}

@article{Puthe2014,
  author  = {P{\"u}the, C. and Kuvshinov, A.},
  title   = {Mapping 3-{D} mantle electrical conductivity from space: a new 3-{D} inversion scheme based on analysis of matrix {\(Q\)}-responses},
  journal = {Geophysical Journal International},
  volume  = {197},
  pages   = {768--784},
  year    = {2014},
  doi     = {10.1093/gji/ggu027}
}

@article{Puthe2015,
  author  = {P{\"u}the, C. and Kuvshinov, A. and Olsen, N.},
  title   = {Handling complex source structures in global {EM} induction studies: from {\(C\)}-responses to new arrays of transfer functions},
  journal = {Geophysical Journal International},
  volume  = {201},
  pages   = {318--328},
  year    = {2015},
  doi     = {10.1093/gji/ggv021}
}

@article{SoleimaniLionheart2006,
  author  = {Soleimani, M. and Lionheart, W. R. B.},
  title   = {Absolute conductivity reconstruction in magnetic induction tomography using a nonlinear method},
  journal = {IEEE Transactions on Medical Imaging},
  volume  = {25},
  pages   = {1521--1530},
  year    = {2006},
  doi     = {10.1109/TMI.2006.884196}
}

@article{SylvesterUhlmann1987,
  author  = {Sylvester, J. and Uhlmann, G.},
  title   = {A global uniqueness theorem for an inverse boundary value problem},
  journal = {Annals of Mathematics},
  volume  = {125},
  pages   = {153--169},
  year    = {1987},
  doi     = {10.2307/1971291}
}

@article{Uhlmann2009,
  author  = {Uhlmann, G.},
  title   = {Electrical impedance tomography and {Calder\'on}'s problem},
  journal = {Inverse Problems},
  volume  = {25},
  pages   = {123011},
  year    = {2009},
  doi     = {10.1088/0266-5611/25/12/123011}
}

@article{Weidelt1972,
  author  = {Weidelt, P.},
  title   = {The inverse problem of geomagnetic induction},
  journal = {Zeitschrift f{\"u}r Geophysik},
  volume  = {38},
  pages   = {257--289},
  year    = {1972}
}

@article{Weidelt2005,
  author  = {Weidelt, P.},
  title   = {The relationship between the spectral function and the underlying conductivity structure in 1-{D} magnetotellurics},
  journal = {Geophysical Journal International},
  volume  = {161},
  pages   = {566--590},
  year    = {2005},
  doi     = {10.1111/j.1365-246X.2005.02625.x}
}

\appendix
\section{Radial conductivity: scalar responses and inverse spectral data}\label{app:radial}
The classical radial problem provides useful context but has a different mathematical structure from the inverse problem for laterally varying conductivity developed in the main text. If $\sigma=\sigma(r)$, the response is diagonal in spherical harmonics, with scalar entry $q_\ell(s)$. For degree $\ell\geq1$, write the normalized electric field $\eg=\bfE/(-\mu s)$ as a toroidal vector spherical harmonic times $w_\ell(r)/r$. Its regular radial factor satisfies
\begin{equation}\label{eq:radialODE}
 -w_\ell''+\left[\frac{\ell(\ell+1)}{r^2}+s\mu\sigma(r)\right]w_\ell=0,
 \qquad w_\ell(r)\sim c\,r^{\ell+1}\quad(r\downarrow0).
\end{equation}
Define the $C$-response, which has dimensions of length, by $C_\ell(s)=w_\ell(R)/w_\ell'(R)$. Set $\lambda_\ell=\ell(\ell+1)$, normalize the incident potential coefficient to one, and choose $\hat x\times\nabla_{\mathbb{S}^2}Y_{\ell m}/\sqrt{\lambda_\ell}$ as the toroidal angular factor. Matching the normal and tangential components of the curl of the normalized electric field with the magnetic Gauss coefficients yields
\begin{equation}
 \frac{\sqrt{\lambda_\ell}\,w_\ell(R)}{R}
 =\ell-(\ell+1)q_\ell(s),\qquad
 \frac{w_\ell'(R)}{\sqrt{\lambda_\ell}}=1+q_\ell(s).
\end{equation}
Consequently,
\begin{equation}\label{eq:QC}
 C_\ell(s)=\frac{R}{\ell+1}\frac{1-\frac{\ell+1}{\ell}q_\ell(s)}{1+q_\ell(s)},
 \qquad
 q_\ell(s)=\frac{\frac{\ell}{\ell+1}-\ell C_\ell(s)/R}{1+\ell C_\ell(s)/R}.
\end{equation}
These ratios may be interpreted meromorphically; $q_\ell(0)=0$ and $C_\ell(0)=R/(\ell+1)$. The formula connects global single-harmonic responses, not arbitrary stationwise transfer functions in the presence of mixed external sources \cite{Puthe2015}.

A change of variables makes the inverse spectral character explicit. Set $v=w_\ell/r^{\ell+1}$,
\begin{equation}\label{eq:radialt}
 t=\frac{R}{2\ell+1}\left[(R/r)^{2\ell+1}-1\right],\qquad
 u(t,s)=v(r(t),s).
\end{equation}
Then
\begin{equation}\label{eq:string}
 -u_{tt}=(-s)b_\ell(t)u,
 \qquad b_\ell(t)=\mu\sigma(r(t))(r(t)/R)^{4\ell+4},
\end{equation}
with the boundary relation
\begin{equation}\label{eq:weyl}
 -\frac{u_t(0,s)}{u(0,s)}=\frac1{C_\ell(s)}-\frac{\ell+1}{R}.
\end{equation}
For bounded $\sigma$, the regular branch has $v'(r)=O(r)$ as $r\downarrow0$, hence
\[
 \int_0^\infty |u_t|^2\dd t
 =\int_0^R |v'(r)|^2(r/R)^{2\ell+2}\dd r<\infty.
\] For nonreal $s$, this identifies it as the Weyl solution \cite[equation~(2.7) and Lemma~4.2]{EckhardtKostenko2016}. This string belongs to the nonnegative class of generalized strings in \cite[Theorem~6.1 and Proposition~7.3]{EckhardtKostenko2016}, with infinite length, mass measure $b_\ell(t)\dd t$, and vanishing quadratic spectral term. In the normalization of equation~(5.1) of that reference,
\begin{equation}\label{eq:string-Weyl-normalization}
 m_\ell(-s)=\frac1s\left[\frac1{C_\ell(s)}-\frac{\ell+1}{R}\right],\qquad s\ne0.
\end{equation}
By the inverse string theorem, $m_\ell$ uniquely determines $b_\ell(t)$. As $t\mapsto r$ is known, this recovers $\sigma(r)$ almost everywhere from the full frequency-dependent $q_\ell$ for any fixed $\ell\geq1$. The inverse spectral viewpoint is part of the classical one-dimensional induction literature \cite{Bailey1970,Weidelt1972,Weidelt2005}.

To derive \eqref{eq:radialslope}, rescale the regular solution so that $v(R,s)=1$ near $s=0$, leaving $C_\ell$ and $q_\ell$ unchanged. Integrating the equation for $v$ yields
\begin{equation}\label{eq:radial-logderivative}
 \frac1{C_\ell(s)}-\frac{\ell+1}{R}
 =v'(R,s)=\frac{s\mu}{R^{2\ell+2}}\int_0^R\sigma(r)r^{2\ell+2}v(r,s)\dd r.
\end{equation}
At $s=0$ one has $v(r,0)=1$. Differentiating \eqref{eq:radial-logderivative} with respect to $s$ at $s=0$ and using \eqref{eq:QC} gives \eqref{eq:radialslope}. The inverse spectral argument uses the full frequency dependence at a fixed degree, whereas the radial and axisymmetric moment arguments use zero-frequency slopes across spatial harmonic degrees.

\section{Technical details for single-frequency uniqueness}\label{app:single-frequency}

This appendix supplies the functional-analytic and asymptotic details used in
Section~\ref{subsec:shell}. All five-dimensional functions and operators are
restricted to the SO(4)-invariant subspace whenever this is indicated. For
variational identities we use the usual sesquilinear pairing, linear in the
first argument. In the DN asymptotic identity in Lemma~\ref{lem:outer}, we
explicitly use the bilinear pairing without complex conjugation.

\subsection{Energy correspondence for the five-dimensional lifting}
\label{app:5d-energy}

Let $U=U(\rho,z)$ be smooth and radial in the four variables
$\mathbf y\in\mathbb{R}^4$, with $\rho=|\mathbf y|$, and define
$\eg=\rho U\ephi$ in $\mathbb{B}_R^3$. In cylindrical coordinates,
\[
 \dd x=\rho\,\dd\rho\,\dd\varphi\,\dd z,
 \qquad
 \dd X=\rho^3\,\dd\rho\,\dd\Omega_3\,\dd z,
 \qquad |\mathbb{S}^3|=2\pi^2.
\]
Consequently,
\[
 \int_{\mathbb{B}_R^3}|\eg|^2\dd x
 =2\pi\int_{D_R^{\mathrm{mer}}}\rho^3|U|^2\dd\rho\dd z
 =\frac1\pi\int_{\mathbb{B}_R^5}|U|^2\dd X,
\]
where $D_R^{\mathrm{mer}}=\{(\rho,z):\rho>0,\ \rho^2+z^2<R^2\}$.
A direct calculation gives
\[
 \curl\eg=-\rho U_z\mathbf e_\rho+(2U+\rho U_\rho)\mathbf e_z.
\]
Hence
\[
 \int_{\mathbb{B}_R^3}|\curl\eg|^2\dd x
 =2\pi\int_{D_R^{\mathrm{mer}}}
 \left(\rho^3|U_z|^2+\rho|2U+\rho U_\rho|^2\right)
 \dd\rho\dd z.
\]
For smooth $U$, integration by parts in $\rho$ yields
\[
 \int_{D_R^{\mathrm{mer}}}
 \left(\rho|2U+\rho U_\rho|^2-\rho^3|U_\rho|^2\right)
 \dd\rho\dd z
 =
 2\int_{-R}^{R}(R^2-z^2)|U|^2_{\rho=\sqrt{R^2-z^2}}\dd z.
\]
The surface measure on $\partial \mathbb{B}_R^5$ in the coordinates
$(\rho,z)$ is
\[
 \dd S_5=2\pi^2R\rho^2\dd z.
\]
Combining these formulas gives
\[
 \int_{\mathbb{B}_R^3}|\curl\eg|^2\dd x
 =
 \frac1\pi\left(
 \int_{\mathbb{B}_R^5}|\nabla U|^2\dd X+
 \frac2R\int_{\partial \mathbb{B}_R^5}|U|^2\dd S_5
 \right).
\]
Polarization gives the corresponding sesquilinear identities for two
invariant functions $U,V$. Smooth invariant functions are dense in the
invariant $H^1(\mathbb{B}_R^5)$ space. The identities therefore extend by density
to the natural energy spaces. In particular, the map
$U\mapsto \rho U\ephi$ is continuous from invariant $H^1(\mathbb{B}_R^5)$ into
the axisymmetric toroidal subspace of $H(\curl;\mathbb{B}_R^3)$.

\subsection{Boundary trace isomorphisms}
\label{app:trace-isomorphisms}

To distinguish unit-sphere normalization from normalization at radius $R$, set
\[
 Y_\ell^{(R)}=R^{-1}Y_{\ell0},\qquad
 \Psi_{\ell-1}^{(R)}=R^{-2}\Psi_{\ell-1},\qquad \ell\geq1,
\]
where the angular arguments are understood. Here $Y_{\ell0}$ is the
real $L^2(\mathbb{S}^2)$-normalized zonal harmonic proportional to
$P_\ell(\xi)$, and $\Psi_{\ell-1}$ is the real
$L^2(\mathbb{S}^4)$-normalized $\mathrm{SO}(4)$-invariant harmonic
proportional to $C_{\ell-1}^{(3/2)}(\xi)$, with positive normalization
constants. Thus $Y_\ell^{(R)}$ and $\Psi_{\ell-1}^{(R)}$ are normalized on
$\partial\mathbb{B}_R^3$ and $\partial\mathbb{B}_R^5$, respectively.
The surface measures in the angular variable $\xi$ are
\[
 \dd S_3=2\pi R^2\dd\xi,\qquad
 \dd S_5=2\pi^2R^4(1-\xi^2)\dd\xi.
\]
The identity
$P_\ell'(\xi)=C_{\ell-1}^{(3/2)}(\xi)$ and the standard Legendre norm
formula imply
\[
 \|P_\ell'\|_{L^2(\partial\mathbb{B}_R^5)}^2
 =
 \pi R^2\ell(\ell+1)
 \|P_\ell\|_{L^2(\partial\mathbb{B}_R^3)}^2.
\]
With the positive normalization constants, this gives
\[
 \mathcal D_RY_\ell^{(R)}
 =-\frac{\sqrt{\pi\ell(\ell+1)}}R\Psi_{\ell-1}^{(R)},
 \qquad
 \mathcal M_R\Psi_{\ell-1}^{(R)}
 =\frac{\sqrt{\ell(\ell+1)}}{\sqrt\pi R}Y_\ell^{(R)}.
\]
The eigenvalues of the negative Laplace--Beltrami operators on the two
spheres are $\ell(\ell+1)/R^2$ and $(\ell-1)(\ell+2)/R^2$, respectively.
Thus the $H^{1/2}$ and $H^{-1/2}$ Sobolev norms are equivalent to the
corresponding weighted coefficient norms. Since the multipliers above
grow like $\ell/R$, they compensate exactly for the change from
$H^{1/2}$ to $H^{-1/2}$. Therefore
\[
 \mathcal D_R:
 H_{\mathrm{zon},0}^{1/2}(\partial\mathbb{B}_R^3)
 \xrightarrow{\ \sim\ }
 H_{\mathrm{SO}(4)}^{-1/2}(\partial\mathbb{B}_R^5),
\]
and
\[
 \mathcal M_R:
 H_{\mathrm{SO}(4)}^{1/2}(\partial\mathbb{B}_R^5)
 \xrightarrow{\ \sim\ }
 H_{\mathrm{zon},0}^{-1/2}(\partial\mathbb{B}_R^3).
\]
The index shift is important: the three-dimensional dipole $\ell=1$
is mapped to the constant $\mathrm{SO}(4)$-invariant harmonic
$\Psi_0^{(R)}$ on $\partial\mathbb{B}_R^5$. Hence no invariant boundary
mode is omitted.

\subsection{Weak realization of the Maxwell boundary traces}
\label{app:weak-maxwell-traces}

Let $h=\mathcal D_R\tau$, and let $U\in H^1(\mathbb{B}_R^5)$ be the SO(4)-invariant
solution of
\[
 (-\Delta_5+q)U=0,\qquad
 \left(\partial_r+\frac2R\right)U=h.
\]
The weak Robin formulation is
\begin{equation}\label{eq:weak-robin-5d}
 \int_{\mathbb{B}_R^5}
 \left(\nabla U\cdot\nabla\overline V+qU\overline V\right)\dd X
 +\frac2R\int_{\partial \mathbb{B}_R^5}u\overline v\dd S_5
 =
 \langle h,v\rangle_{\partial \mathbb{B}_R^5}
\end{equation}
for every SO(4)-invariant $V\in H^1(\mathbb{B}_R^5)$, where
$v=V|_{\partial \mathbb{B}_R^5}$ and $u=U|_{\partial \mathbb{B}_R^5}$.

Every smooth axisymmetric toroidal vector field in $\mathbb{B}_R^3$ has the form
$\mathbf v=v_\varphi(\rho,z)\mathbf e_\varphi=\rho V(\rho,z)\mathbf e_\varphi$,
where regularity across the axis implies that $V$ extends as a smooth even
function of the signed radial variable near $\rho=0$. Equivalently, locally
$V(\rho,z)=F(\rho^2,z)$ with $F$ smooth. Hence
$V_5(\mathbf y,z)=F(|\mathbf y|^2,z)$ is a smooth SO(4)-invariant scalar
function. Conversely, every smooth SO(4)-invariant scalar test function
produces a smooth axisymmetric toroidal vector field by this lifting.
Thus the scalar tests used below exhaust the smooth axisymmetric toroidal
test fields, rather than merely providing a subclass of them.

For a smooth invariant $V$, define the three-dimensional azimuthal
test field $\mathbf v=\rho V\ephi$. By the polarized form of the
energy identities in Appendix~\ref{app:5d-energy},
\[
 \begin{aligned}
 &\int_{\mathbb{B}_R^3}
 \left(\curl\eg\cdot\curl\overline{\mathbf v}
       +s\mu\sigma\eg\cdot\overline{\mathbf v}\right)\dd x\\
 &\qquad=
 \frac1\pi\left[
 \int_{\mathbb{B}_R^5}
 \left(\nabla U\cdot\nabla\overline V+qU\overline V\right)\dd X
 +\frac2R\int_{\partial \mathbb{B}_R^5}u\overline v\dd S_5
 \right]\\
 &\qquad=\frac1\pi
 \langle h,v\rangle_{\partial \mathbb{B}_R^5}.
 \end{aligned}
\]
Since $\bfH=\curl\eg$ and
$\curl\bfH+s\mu\sigma\eg=0$, the curl Green formula gives
\[
 \int_{\mathbb{B}_R^3}
 \left(\curl\eg\cdot\curl\overline{\mathbf v}
       +s\mu\sigma\eg\cdot\overline{\mathbf v}\right)\dd x
 =
 \int_{\partial \mathbb{B}_R^3}
 (\bfH\times\nu)\cdot\overline{\mathbf v_T}\dd S_3.
\]
On the other hand, using $h=-R^{-2}\partial_\xi\tau$ and the two
surface measures,
\[
 \frac1\pi\langle h,v\rangle_{\partial \mathbb{B}_R^5}
 =
 -2\pi R^2\int_{-1}^{1}
 (1-\xi^2)\partial_\xi\tau\,\overline V\dd\xi.
\]
This is exactly
\[
 \int_{\partial \mathbb{B}_R^3}
 \bigl[(-\nabla_\Gamma\tau)\times\nu\bigr]
 \cdot\overline{\mathbf v_T}\dd S_3.
\]
Thus the weak tangential magnetic trace is
$\bfH_T=-\nabla_\Gamma\tau$ for all azimuthal tests.

To see that this verifies the full weak boundary condition, average an
arbitrary smooth vector test field over rotations about the $z$ axis.
Because the coefficient and the constructed field are axisymmetric, the
weak pairing is unchanged by this averaging; equivalently, all nonzero
azimuthal Fourier modes pair to zero. The averaged test field splits into
an azimuthal component and a meridional component. The meridional component
has zero pairing: its curl is azimuthal, whereas $\curl\eg$ is
meridional, and it is pointwise orthogonal to $\eg$; the boundary pairing
with $\bfH\times\nu$ also vanishes. The azimuthal component is of the
form established above and is therefore represented by a smooth invariant
scalar test. Hence the scalar tests suffice for the full weak formulation.

Finally, for a smooth zonal scalar test function $g$,
\[
 \begin{aligned}
 \langle H_\nu,g\rangle_{\partial \mathbb{B}_R^3}
 &=\int_{\partial \mathbb{B}_R^3}
   (\nu\times\eg)\cdot\nabla_\Gamma\overline g\dd S_3\\
 &=2\pi R^2\int_{-1}^{1}
   (1-\xi^2)u\,\partial_\xi\overline g\dd\xi\\
 &=\langle\mathcal M_Ru,g\rangle_{\partial \mathbb{B}_R^3}.
 \end{aligned}
\]
The identity extends to the stated Sobolev trace spaces by density.
This verifies both the tangential input and the normal magnetic output
used in \eqref{eq:Robin-Dirichlet-map}.

\subsection{Compact-resolvent estimate for high-degree solutions}
\label{app:compact-resolvent}

Let $q\in L^\infty(\mathbb{B}_1^5)$ be purely imaginary and let
$\mathcal{A}_q=-\Delta_{5,D}+q$. The Dirichlet variational form is coercive
because
\[
 \operatorname{Re}\int_{\mathbb{B}_1^5}
 \left(|\nabla w|^2+q|w|^2\right)\dd X
 =\|\nabla w\|_2^2.
\]
Thus $\mathcal{A}_q^{-1}:L^2(\mathbb{B}_1^5)\to H_0^1(\mathbb{B}_1^5)$ is bounded. Since
$H_0^1(\mathbb{B}_1^5)\hookrightarrow L^2(\mathbb{B}_1^5)$ is compact,
$\mathcal{A}_q^{-1}$ is compact as an operator on $L^2$.

The functions $v_n=r^n\Psi_n$ are mutually orthogonal in $L^2(\mathbb{B}_1^5)$
and satisfy $\|v_n\|_2=(2n+5)^{-1/2}$. Hence
$w_n=v_n/\|v_n\|_2$ converges weakly to zero. Multiplication by $q$
is bounded on $L^2$, so $qw_n\rightharpoonup0$. Compactness gives
\[
 \|\mathcal{A}_q^{-1}(qw_n)\|_2\longrightarrow0.
\]
Using
\[
 u_n-v_n=-\mathcal{A}_q^{-1}(qv_n)
        =-\|v_n\|_2\mathcal{A}_q^{-1}(qw_n),
\]
we obtain
\[
 \|u_n-v_n\|_2=o(\|v_n\|_2)=o(n^{-1/2}),
\]
which is the estimate used in Lemma~\ref{lem:outer}.

\subsection{The Gegenbauer Jacobi-matrix limit}
\label{app:gegenbauer-limit}

Let $\{\Psi_n\}_{n\geq0}$ be the real orthonormal SO(4)-invariant
spherical harmonics on $\mathbb{S}^4$. The normalized Gegenbauer
three-term recurrence is
\[
 \xi\Psi_n=a_n\Psi_{n+1}+a_{n-1}\Psi_{n-1},
 \qquad
 a_n=\sqrt{\frac{(n+1)(n+3)}{(2n+3)(2n+5)}}.
\]
For a polynomial $p$, the matrix element
$\langle p(\xi)\Psi_n,\Psi_{n+d}\rangle$ is a finite linear
combination of products of finitely many coefficients
$a_{n+j}$. For fixed $d$ and fixed polynomial degree, all shifts
$j$ remain bounded as $n\to\infty$. Since $a_{n+j}\to1/2$,
the limit is the $(d,0)$ matrix coefficient of the constant
Jacobi operator on $\ell^2(\mathbb Z)$ with off-diagonal entries
$1/2$.

Under the Fourier transform on $\ell^2(\mathbb Z)$, this operator
is multiplication by $\cos\theta$. Therefore
\[
 \lim_{n\to\infty}
 \langle p(\xi)\Psi_n,\Psi_{n+d}\rangle_{L^2(\mathbb{S}^4)}
 =
 \frac1\pi\int_0^\pi p(\cos\theta)\cos(d\theta)\dd\theta.
\]
For general $p\in C([-1,1])$, choose polynomials $p_m$ with
$\|p-p_m\|_\infty\to0$. Uniformly in $n$,
\[
 \left|
 \langle(p-p_m)\Psi_n,\Psi_{n+d}\rangle
 \right|
 \leq\|p-p_m\|_\infty.
\]
The same bound holds for the limiting integral. Passing first
$n\to\infty$ for fixed $m$ and then $m\to\infty$ proves the
continuous-function version used in Lemma~\ref{lem:outer}.

\subsection{Density of reachable interface traces}
\label{app:layer-stripping}

We give the functional-analytic proof of the density property
\eqref{eq:reachable-trace-density}. Let
$\Omega=\mathbb{B}_R^5$, $\Sigma=\partial \mathbb{B}_a^5$, and let
$q\in L^\infty(\Omega)$ be purely imaginary and SO(4)-invariant.
Define the interior trace map
\[
 \mathcal E_q:H_{\mathrm{SO}(4)}^{1/2}(\partial\Omega)
 \longrightarrow H_{\mathrm{SO}(4)}^{1/2}(\Sigma),
 \qquad
 \mathcal E_qf=u_f|_\Sigma,
\]
where $u_f$ is the unique global solution of
$(-\Delta_5+q)u_f=0$ with outer Dirichlet value $f$.
The trace theorem and the bounded solution operator show that
$\mathcal E_q$ is bounded. We prove that its range is dense.

Suppose $\eta\in H_{\mathrm{SO}(4)}^{-1/2}(\Sigma)$ satisfies
\[
 \langle\eta,\mathcal E_qf\rangle_\Sigma=0
 \qquad\text{for all }f\in
 H_{\mathrm{SO}(4)}^{1/2}(\partial\Omega),
\]
where the duality is linear in the first argument and conjugate-linear
in the second. Let $\mathcal{P}_{\mathrm{SO}(4)}$ denote averaging over the SO(4)
action on the $\mathbf y$ variables. We interpret the interface source
$\eta\delta_\Sigma$ on the full space $H_0^1(\Omega)$ by
\[
 \phi\longmapsto
 \langle\eta,\mathcal{P}_{\mathrm{SO}(4)}(\phi|_\Sigma)\rangle_\Sigma.
\]
This is a continuous functional by the trace theorem and boundedness
of the averaging projection.

Solve the adjoint interface-source problem
\[
 (-\Delta_5+\overline q)v=\eta\delta_\Sigma
 \quad\text{in }\Omega,
 \qquad v|_{\partial\Omega}=0.
\]
More precisely, $v\in H_0^1(\Omega)$ is defined by
\begin{equation}\label{eq:reachable-adjoint-weak}
 \int_\Omega
 \left(\nabla v\cdot\nabla\overline\phi
       +\overline q\,v\overline\phi\right)\dd X
 =\langle\eta,\mathcal{P}_{\mathrm{SO}(4)}(\phi|_\Sigma)\rangle_\Sigma
 \qquad(\phi\in H_0^1(\Omega)).
\end{equation}
The real part of the form at $(v,v)$ is $\|\nabla v\|_2^2$,
so Poincar\'e's inequality and Lax--Milgram give a unique solution.
Invariance of $q$ and of the source, together with uniqueness, implies
that $v$ is SO(4)-invariant. Consequently the projection
$\mathcal{P}_{\mathrm{SO}(4)}$ can be omitted when the test function is invariant;
by averaging, the weak equation also holds for arbitrary test functions.

For any global solution $u_f$, apply the second Green identity to the
adjoint solution $v$ and $u_f$ on the whole ball, with the interface source
interpreted distributionally. The outward normal $\nu$ on
$\partial\Omega$ is the normal of the full ball. Since
$(-\Delta_5+q)u_f=0$, $v=0$ on $\partial\Omega$, and the interface
source is $\eta\delta_\Sigma$, the identity gives
\begin{equation}\label{eq:reachable-green}
 \langle\eta,u_f|_\Sigma\rangle_\Sigma
 =-\langle\partial_\nu v|_{\partial\Omega},f\rangle_{\partial\Omega}.
\end{equation}
The assumed annihilation therefore implies that the outer Neumann trace
of $v$ vanishes in the invariant trace space. Since $v$ is invariant,
this is its full outer Neumann trace. Its outer Dirichlet trace is zero
by construction.

In the shell $\Omega\setminus\overline{\mathbb{B}_a^5}$, the interface source
is absent and $v$ solves the homogeneous equation with bounded
potential $\overline q$. Since both the Dirichlet and Neumann traces of $v$ vanish on
$\partial\Omega$, extension by zero across the outer sphere is a
weak solution with a bounded potential extended by zero. The unique-continuation
theorem for Schr\"odinger operators with bounded potentials \cite{JerisonKenig1985}
gives $v=0$ in the shell. Hence $v|_\Sigma=0$.
Inside $\mathbb{B}_a^5$, the function $v$ solves the homogeneous adjoint
Dirichlet problem with zero boundary trace. Coercivity gives $v=0$
there as well. Thus $v$ vanishes on both sides of $\Sigma$ and is zero
in $\Omega$. Equation~\eqref{eq:reachable-adjoint-weak} then shows
that $\eta$ annihilates the traces on $\Sigma$ of all functions in
$H_0^1(\Omega)$. The interior-interface trace map
$H_0^1(\Omega)\to H^{1/2}(\Sigma)$ is surjective and has a bounded
right inverse, since $\Sigma$ is a smooth closed hypersurface strictly
inside $\Omega$. The same property holds on the SO(4)-invariant
subspaces by composing a right inverse with the averaging projection
$\mathcal{P}_{\mathrm{SO}(4)}$. Consequently $\eta=0$ in
$H_{\mathrm{SO}(4)}^{-1/2}(\Sigma)$.

We have shown that the annihilator of $\operatorname{Ran}\mathcal E_q$
is trivial. By the dual characterization of dense subspaces,
\[
 \overline{\operatorname{Ran}\mathcal E_q}
 =H_{\mathrm{SO}(4)}^{1/2}(\Sigma),
\]
which proves the density property \eqref{eq:reachable-trace-density}.

\end{document}